\documentclass[sn-mathphys,Numbered]{sn-jnl}
\pdfoutput=1
\usepackage{pgfplots}
\usepackage{graphicx}%
\usepackage{multirow}%
\usepackage{amsmath,amssymb,amsfonts}%
\usepackage{amsthm}%
\usepackage{mathrsfs}%
\usepackage[title]{appendix}%
\usepackage{xcolor}%
\usepackage{textcomp}%
\usepackage{manyfoot}%
\usepackage{booktabs}%
\usepackage{algpseudocode}%
\usepackage{listings}%
\usepackage{ragged2e}
\usepackage{fullpage}
\usepackage{changepage}
\usepackage{stfloats}
\usepackage[caption=false]{subfig}
\usepackage{extarrows}
\usepackage{cleveref}
\usepackage{tikz}
\usepackage{pifont}
\usepackage{enumitem}
\usepackage[ruled, vlined, noline]{algorithm2e}
\usepackage{float}
\usepackage{graphicx}
\usepackage{subcaption}
\theoremstyle{thmstyleone}%
\newtheorem{theorem}{Theorem}[section]%  meant for continuous numbers
\newtheorem{proposition}[theorem]{Proposition}% 
\theoremstyle{thmstyletwo}%
\newtheorem{remark}{Remark}%
\theoremstyle{thmstylethree}%
\newtheorem{lemma}{Lemma}[section]
\AtBeginDocument{
	\setlength{\abovedisplayskip}{6pt}
	\setlength{\belowdisplayskip}{6pt}
	\setlength{\abovedisplayshortskip}{3pt}
	\setlength{\belowdisplayshortskip}{3pt}
	
	\setlength{\jot}{2pt}
}
\begin{document}
	\title[Article Title]{Convergence rates of mixed primal-dual dynamical systems with implicit Hessian damping and Tikhonov regularization}
	
	\author[1]{\fnm{Hong-lu} \sur{Li}}\email{lihonglu$\_$7988@163.com}
	
	\author*[1]{\fnm{Yi-bin} \sur{Xiao}}\email{xiaoyb9999@hotmail.com}
	
	\author[2]{\fnm{Xin} \sur{He}}\email{hexinuser@163.com}
	
	\author[3]{\fnm{Kequan} \sur{Zhao}}\email{kequanz@163.com}
	
	\affil[1]{\orgdiv{School of Mathematical Sciences}, \orgname{University of Electronic Science and Technology of China}, \orgaddress{\street{}\city{Chengdu}, \postcode{611731}, \state{Sichuan}, \country{P.R. China}}}
	
	\affil[2]{\orgdiv{School of Science}, \orgname{Xihua University}, \orgaddress{\street{}\city{Chengdu}, \postcode{610039}, \state{Sichuan}, \country{P.R. China}}}
	
	\affil[3]{\orgdiv{School of Mathematical Sciences}, \orgname{Chongqing Normal University}, \orgaddress{\street{}\city{Chongqing}, \postcode{401331}, \state{} \country{P.R. China}}}
	
	\maketitle

\begin{abstract}
	
	\hspace{2em}	\textbf{Abstract.} We propose a mixed-order primal--dual dynamical system with implicit Hessian damping and Tikhonov regularization for linearly constrained convex optimization, and derive a corresponding fast primal--dual algorithm. The proposed dynamical system couples a second-order primal equation with a first-order dual equation without requiring twice differentiability of the objective function. An energy analysis establishes $\mathcal{O}(t^{-2})$ convergence rates for the augmented Lagrangian gap, feasibility violation, and objective residual along an extrapolated primal trajectory. For the proposed algorithm, we construct discrete energy estimates and establish the corresponding $\mathcal{O}(k^{-2})$ rates at the extrapolated primal iterates. In particular, in the critical regularization regimes
	$\epsilon(t)=a/t^2$ and $\epsilon_k=a/k^2$, these accelerated rates hold simultaneously with strong convergence of the primal trajectory
	and primal iterates to the minimum-norm solution, respectively, under suitable parameter and trajectory conditions. Numerical experiments illustrate the convergence behavior of the proposed dynamical system and the corresponding algorithm.
	
\end{abstract}

%All acknowledgements should be placed in the back of the paper after Conclusions..

\section{Introduction}\label{sec1}
We consider the linearly constrained convex optimization problem
\begin{equation}\label{problem1}
	\min_{x\in\mathbb{R}^n} f(x)
	\quad \text{subject to} \quad Ax=b,
\end{equation}
where $f:\mathbb{R}^n\to\mathbb{R}$ is convex and continuously
differentiable with a Lipschitz continuous gradient,
$A:\mathbb{R}^n\to\mathbb{R}^m$ is a linear operator, and
$b\in\mathbb{R}^m$.
We denote the feasible set by
$\mathbb{F}:=\{x\in\mathbb{R}^n:Ax=b\}$.
Problems of this form arise in sparse recovery, image reconstruction,
distributed optimization, and PDE-constrained optimization
\cite{Becker SIAM JIS 2011,Herman RTI 1995,
	Nedic IEEE TAC 2014,Rees SIAM JSC 2010}.
When $A=0$ and $b=0$, problem~\eqref{problem1} reduces to
\begin{equation}\label{unconstrained_problem}
	\min_{x\in\mathbb{R}^n} f(x).
\end{equation}

Dynamical systems provide a useful framework for understanding
and designing accelerated optimization algorithms.
For problem~\eqref{unconstrained_problem},
Su et al.~\cite{Su JMLR 2016} studied the inertial dynamical system
with asymptotically vanishing damping
\[
(\mathrm{AVD}_{\alpha}) \qquad
\ddot{x}(t)+\frac{\alpha}{t}\dot{x}(t)+\nabla f(x(t))=0.
\]
For $\alpha=3$, this dynamical system arises as a continuous-time
limit of Nesterov's accelerated gradient method
\cite{Nesterov IMD 1983}.
Under suitable assumptions, its trajectories satisfy
$f(x(t))-\min f=\mathcal{O}(t^{-2})$, reflecting the accelerated
objective decrease of the discrete method.
Subsequent studies developed energy-based and variational
approaches to accelerated dynamical systems and algorithms,
establishing convergence rates and trajectory or iterate
convergence under various assumptions; see
\cite{Attouch MP 2018,Attouch SIAM JO 2016,
	Attouch SIAM JO 2019,Wibisono PNAS 2016,Wilson JMLR 2021}.

For problem~\eqref{problem1}, the Lagrangian framework leads to
primal--dual dynamical systems that jointly evolve the primal
and dual variables.
Zeng et al.~\cite{Zeng IEEE TAC 2023} proposed a second-order
primal--dual dynamical system and established, under appropriate
parameter conditions,
$f(x(t))-f(x^\ast)=\mathcal{O}(t^{-2})$ and
$\|Ax(t)-b\|^2=\mathcal{O}(t^{-2})$.
Bo\c{t} and Nguyen~\cite{Bot JDE 2021} studied an
augmented-Lagrangian-based second-order primal--dual dynamical
system with vanishing damping, obtaining fast convergence rates
and trajectory convergence.
Building on this dynamical system,
Bo\c{t} et al.~\cite{Bot MP 2023} derived a fast augmented
Lagrangian method with $\mathcal{O}(k^{-2})$ convergence rates
and convergence guarantees for the iterates.
He et al.~\cite{He NA 2022} also developed inertial accelerated
primal--dual methods with $\mathcal{O}(k^{-2})$ estimates for
the objective residual and feasibility violation.
Related accelerated convergence results based on inertia,
damping, and time scaling were established in
\cite{He SIAM JCO 2021,He AA 2023,Attouch JOTA 2022,
	Bitterlich NFAO 2020,Bot JDE 2020}.
Other studies obtained exponential convergence or stability
results under additional structural assumptions
\cite{Luo ESAIM COCV 2022,Tang SCL 2020,
	Liang Automatica 2019,Hassan Automatica 2021}. In the second-order primal--dual dynamical systems studied in
\cite{Zeng IEEE TAC 2023,Bot JDE 2021}, both the primal and dual
variables satisfy second-order differential equations, which may
increase discretization cost.
This motivates considering mixed-order primal--dual dynamical
systems that retain inertia in the primal equation while using
a first-order dual equation.
He et al.~\cite{He IEEE TAC 2022} studied such a dynamical system
with time scaling and established fast convergence rates
without strong convexity.
He et al.~\cite{He Automatica 2022} subsequently derived a
primal--dual algorithm through discretization.
Under suitable scaling and parameter conditions, the algorithm
achieves accelerated decay of both the objective residual and
feasibility violation, including an $\mathcal{O}(k^{-2})$ rate.

However, inertial dynamical systems with vanishing damping may
exhibit pronounced trajectory oscillations.
To mitigate this behavior, Hessian-driven damping introduces
additional dissipation through the term
$\nabla^2 f(x(t))\dot{x}(t)$.
Attouch et al.~\cite{Attouch JDE 2016} studied the dynamical system
\[
\ddot{x}(t)+\frac{\alpha}{t}\dot{x}(t)
+\beta\nabla^2 f(x(t))\dot{x}(t)+\nabla f(x(t))=0.
\]
Under the assumptions of their analysis, they established
$f(x(t))-\min f=\mathcal{O}(t^{-2})$ for $\alpha\geq 3$;
for $\alpha>3$, they obtained the sharper rate $o(t^{-2})$
and weak convergence of the trajectories to a minimizer.
Using the identity
$\nabla^2 f(x(t))\dot{x}(t)=\frac{d}{dt}\nabla f(x(t))$,
Attouch et al.~\cite{Attouch MP 2022} derived inertial
first-order algorithms with $\mathcal{O}(k^{-2})$ objective
residual estimates and quantitative gradient decay estimates.
Shi et al.~\cite{Shi MP 2022} used high-resolution differential
equations to identify gradient-correction terms that help
explain the behavior of accelerated methods, and developed
accelerated algorithms through symplectic discretizations
in \cite{Shi NeurIPS 2019}. An alternative is to incorporate the Hessian damping effect
implicitly by evaluating the gradient at an extrapolated point.
This avoids explicit Hessian evaluations and does not require
twice differentiability of the objective function.
Alecsa et al.~\cite{Alecsa AMOR 2021} studied an inertial
dynamical system with asymptotically vanishing damping,
where the gradient is evaluated at the extrapolated point
$x(t)+(\gamma+\beta/t)\dot{x}(t)$.
A formal first-order Taylor expansion of the gradient
explains the implicit Hessian damping effect.
For convex objectives with Lipschitz continuous gradients,
they established $\mathcal{O}(t^{-2})$ objective residual
estimates at both $x(t)$ and
$x(t)+(\gamma+\beta/t)\dot{x}(t)$, together with trajectory
convergence under suitable parameter conditions.
Their numerical results also illustrated the attenuation
of trajectory oscillations.
Further convergence and stabilization results for
Hessian-driven dynamical systems can be found in
\cite{Alvarez JMPA 2002,Attouch DEA 2012,Chen arXiv 2019,
	Attouch JEMS 2023,Attouch EECT 2023}.

For linearly constrained convex optimization,
Niederl{\"a}nder~\cite{Niederlaender SIAM JCO 2021} proposed
a second-order primal--dual dynamical system with Hessian-driven
damping and established accelerated convergence, as well as
exponential-type convergence under additional curvature
assumptions.
He et al.~\cite{He Optim 2025} incorporated Hessian-driven
damping and time scaling into a mixed-order primal--dual
dynamical system.
They obtained accelerated decay of the objective residual
and feasibility violation and established corresponding
convergence properties for an implicit discretization.
Li et al.~\cite{Li arXiv 2025} studied a fully second-order
primal--dual dynamical system with implicit Hessian damping
and Tikhonov regularization, establishing fast convergence
rates and minimum-norm strong convergence under different
regularization regimes, but did not investigate its
discretization.

Although Hessian-driven damping can attenuate oscillations,
it does not by itself ensure convergence to the minimum-norm
solution when the solution set is not a singleton.
This motivates introducing a vanishing Tikhonov term
$\epsilon(t)x(t)$ to promote minimum-norm selection.
For the unconstrained problem~\eqref{unconstrained_problem},
Attouch et al.~\cite{Attouch JMAA 2018} studied the regularized
inertial dynamical system
\[
\ddot{x}(t)+\frac{\alpha}{t}\dot{x}(t)
+\nabla f(x(t))+\epsilon(t)x(t)=0.
\]
Their analysis revealed the role of the decay rate of
$\epsilon(t)$: sufficiently fast decay preserves the
$\mathcal{O}(t^{-2})$ objective convergence rate, whereas
slower decay enables strong convergence to the minimum-norm
minimizer under additional assumptions.
Bo\c{t} et al.~\cite{Bot MP 2021} combined Hessian-driven
damping with Tikhonov regularization and established fast
objective decay and minimum-norm strong convergence in
appropriate regularization regimes.
Related developments include minimum-norm convergence for
Tikhonov-regularized perturbed heavy-ball dynamical systems
\cite{Alecsa SIAM JO 2021}, strong convergence results for
dynamical systems with asymptotically vanishing damping
\cite{Laszlo JDE 2023,Attouch JDE 2022}, and discrete methods
combining fast objective decrease with minimum-norm selection
\cite{Attouch MMOR 2024,Attouch AMO 2023}.

For problem~\eqref{problem1},
Zhu et al.~\cite{Zhu JCAM 2025} studied a Tikhonov-regularized
second-order primal--dual dynamical system with time scaling,
establishing fast convergence rates for the primal--dual gap,
feasibility violation, and objective residual, as well as
minimum-norm strong convergence under suitable conditions.
Further results on regularized primal--dual dynamical systems
were obtained in \cite{Li JOTA 2025,Li arXiv 2025}.
On the algorithmic side,
Bo\c{t} et al.~\cite{Bot NM 2026} developed a fast
forward--backward method that yields a fully splitting
primal--dual algorithm for linearly constrained convex
optimization, with convergence of the iterates and
last-iterate $o(k^{-1})$ rates for the objective residual
and feasibility violation. Related results on strong convergence and residual decay
for Tikhonov-regularized dynamical systems and their
discrete counterparts in monotone operator settings
were established in
\cite{Bot JComplexity 2026,Bot SIOPT 2026}.

Motivated by the above studies, we aim to develop a mixed-order
primal--dual dynamical system and a corresponding discrete
algorithm that combine accelerated convergence with
minimum-norm selection, while incorporating implicit Hessian
damping to attenuate oscillations.
To this end, we consider the augmented Lagrangian function
associated with problem~\eqref{problem1},
\begin{equation*}
	\mathcal{L}_{\rho}(x,\lambda)
	:= f(x)+\langle \lambda, Ax-b\rangle+\frac{\rho}{2}\|Ax-b\|^{2}
\end{equation*}
with $\rho>0$ and propose the following mixed primal--dual dynamical system with implicit Hessian damping and Tikhonov regularization, referred to as the MPDDS-IHTR for short:
\begin{equation}\label{DS}
	\begin{cases}
		\ddot{x}(t) + \dfrac{\alpha}{t}\dot{x}(t)
		= -\nabla_x\mathcal{L}_{\rho}\big(\hat{x}(t),\lambda(t)\big) - \epsilon(t)\,x(t),\\[0.8ex]
		\dot{\lambda}(t)
		= \dfrac{t}{2\alpha}\Big(\big(2\alpha-3t\beta(t)\epsilon(t)\big)Ax(t)+\eta(t)A\dot{x}(t)
		-3t\beta(t)A\nabla_x\mathcal{L}_{\rho}\big(\hat{x}(t), \lambda(t)\big)-2\alpha b\Big),
	\end{cases}
\end{equation}
where $\hat{x}(t):=x(t)+\beta(t)\dot{x}(t)$ and $\eta(t):=3t\big(1+\dot{\beta}(t)\big)-\alpha\beta(t)$.
Here $\alpha>0$ and $t\ge t_{0}>0$. The function $\beta(\cdot)$ is defined by
$\beta(t):=\gamma+\dfrac{\beta}{t}$ with either $\gamma>0$ and $\beta\in\mathbb{R}$, or $\gamma=0$ and $\beta> 0$.
Moreover, $\xi(\cdot):[t_{0},+\infty)\to(0,+\infty)$ is a continuous function, and
$\epsilon(\cdot):[t_{0},+\infty)\to\mathbb{R}_{+}$ is a $\mathcal{C}^{1}$ non-increasing function satisfying
$\lim_{t\to+\infty}\epsilon(t)=0$. Throughout the paper, we assume that the saddle-point set $\Omega$ of the augmented Lagrangian $\mathcal{L}_{\rho}$ is nonempty.

The main contributions of this work are summarized as follows.
\begin{itemize}
	\item[(i)]
	We propose a mixed-order primal--dual dynamical system
	incorporating implicit Hessian damping and Tikhonov
	regularization for problem~\eqref{problem1}.
	The proposed dynamical system retains a first-order dual
	equation and does not require twice differentiability
	of the objective function.
	By constructing suitable energy functions, we establish
	$\mathcal{O}(t^{-2})$ convergence rates for the augmented
	Lagrangian gap, feasibility violation, and objective
	residual evaluated along
	$\hat{x}(t)=x(t)+\beta(t)\dot{x}(t)$.
	We also establish strong convergence of the primal
	trajectory to the minimum-norm solution under the
	stated assumptions.
	
	\item[(ii)]
	We derive an implicit primal--dual algorithm by
	discretizing the proposed dynamical system and provide
	a discrete convergence analysis.
	In particular, we construct discrete energy estimates
	to establish $\mathcal{O}(k^{-2})$ rates for the augmented
	Lagrangian gap, feasibility violation, and objective
	residual at the extrapolated iterates
	$z_k=x_k+\beta_k(x_k-x_{k-1})$.
	We further establish minimum-norm strong convergence
	of the primal iterates under the stated assumptions.
	Numerical experiments illustrate the convergence
	behavior of the proposed algorithm.
	
	\item[(iii)]
	We analyze the critical regularization regimes
	$\epsilon(t)=a/t^2$ and $\epsilon_k=a/k^2$.
	Under suitable parameter and trajectory conditions,
	the proposed dynamical system and algorithm
	simultaneously achieve accelerated convergence rates
	and strong convergence to the minimum-norm solution.
	These regimes are not covered by the weighted
	integrability and summability conditions
	$\int_{t_0}^{+\infty}t\epsilon(t)\,dt<+\infty$
	and $\sum_{k=1}^{+\infty}k\epsilon_k<+\infty$
	used in the respective rate analyses, and therefore
	require separate energy estimates.
\end{itemize}

The remainder of this paper is organized as follows. In Section \ref{DS_section}, we investigate the convergence rate of the dynamical system~\eqref{DS} and prove that the trajectory generated by the dynamical system \eqref{DS} converges to the minimum-norm solution of problem~\eqref{problem1}. In Section \ref{Algorithm_section}, by applying an appropriate discretization scheme to the dynamical system \eqref{DS}, we derive a fast primal--dual algorithm and establish convergence properties that are consistent with those of the dynamical system~\eqref{DS}. In Section \ref{Numerical}, we present a series of numerical experiments to validate the theoretical results. Finally, Section \ref{Conclusion} concludes the paper with a brief summary of our contributions.

\section{Convergence Analysis of the Proposed Dynamical System}\label{DS_section}
In this section, we establish convergence rates and minimum-norm
strong convergence for the dynamical system~\eqref{DS}.
For prescribed initial data at $t_0>0$, global existence and
uniqueness follow from the Picard--Lindel\"of theorem and
continuation, since the associated first-order vector field
is globally Lipschitz in the state variables, uniformly on
each compact time interval, and has at most linear growth
\cite{Teschl AMS 2012}.
We begin with an energy estimate used throughout the analysis.
\begin{lemma}\label{DS_lem1}
	Let $(x(t), \lambda(t))$ be a solution of the dynamical system $\rm{\eqref{DS}}$ and $(x^*, \kappa)\in\mathbb{F}\times\mathbb{R}^m$. Denote
	\begin{equation}\label{DS_lem1_1}
		\begin{aligned}
			\mathcal{E}_{\kappa}(t):=&t^2(\mathcal{L}_{\rho}(x(t)+\beta(t)\dot{x}(t),\kappa)-\mathcal{L}_{\rho}(x^*, \kappa)+\frac{\epsilon(t)}{2}\Vert x(t)\Vert^2)\\
			&+\frac{1}{2}\Vert \frac{2\alpha}{3}(x(t)-x^*)+t\dot{x}(t)\Vert^2+\frac{\alpha(\alpha-3)}{9}\Vert x(t)-x^*\Vert^2
			+\frac{\alpha}{3}\Vert\lambda(t)-\kappa\Vert^2.
		\end{aligned}
	\end{equation}
	Assume that $t_0>0$, $\alpha>3$, and $\epsilon$ is a continuously differentiable and nonincreasing function  that satisfies the following condition: there exist $M>1$ and $t_1\ge t_0$ such that
	\begin{equation}\label{parameter_condition_1}
		2\dot{\epsilon}(t)+M|(\gamma+\frac{\beta}{t})|\epsilon(t )^2\leq 0,\;(\forall t\ge t_1).
	\end{equation}
	Then, there exist $t_2\ge t_0$  and $r_1>0$ such that for all $t\ge t_2$,
	\begin{equation}\label{DS_lem1_6}
		\begin{aligned}
			\dot{\mathcal{E}}_{\kappa}(t)\leq &\frac{2(3-\alpha)}{3}t(\mathcal{L}_{\rho}(x(t)+\beta(t)\dot{x}(t),\kappa)-\mathcal{L}_{\rho}(x^*,\kappa))
			+\frac{\alpha}{3}t\epsilon(t)\Vert x^*\Vert^2\\
			&-\frac{\alpha}{3}t\epsilon(t)\Vert x(t)-x^*\Vert^2+\frac{3-\alpha}{3}t\epsilon(t)\Vert x(t) \Vert^2-r_1t\Vert \dot{x}(t)\Vert^2.
		\end{aligned}
	\end{equation}
\end{lemma}
\begin{proof}
	Differentiating \eqref{DS_lem1_1}, substituting the dynamical
	system~\eqref{DS}, and using $Ax^*=b$, we obtain
	\begin{equation}\label{DS_lem1_2}
		\begin{aligned}
			\dot{\mathcal{E}}_{\kappa}(t)
			={}&2t\left(
			\mathcal{L}_{\rho}(x(t)+\beta(t)\dot{x}(t),\kappa)
			-\mathcal{L}_{\rho}(x^*,\kappa)
			+\frac{\epsilon(t)}{2}\Vert x(t)\Vert^2
			\right)\\
			&+t^2\left(\dot{\beta}(t)-\frac{\alpha}{t}\beta(t)\right)
			\left\langle
			\nabla_x\mathcal{L}_{\rho}
			(x(t)+\beta(t)\dot{x}(t),\lambda(t)),\dot{x}(t)
			\right\rangle\\
			&-t^2\beta(t)
			\left\Vert\nabla_x\mathcal{L}_{\rho}
			(x(t)+\beta(t)\dot{x}(t),\lambda(t))\right\Vert^2
			+\frac{t^2\dot{\epsilon}(t)}{2}\Vert x(t)\Vert^2\\
			&+\frac{3-\alpha}{3}t\Vert\dot{x}(t)\Vert^2
			-\frac{2\alpha}{3}t
			\left\langle
			\nabla_x\mathcal{L}_{\rho}
			(x(t)+\beta(t)\dot{x}(t),\kappa)
			+\epsilon(t)x(t),x(t)-x^*
			\right\rangle\\
			&-t^2\beta(t)\epsilon(t)
			\left\langle
			\nabla_x\mathcal{L}_{\rho}
			(x(t)+\beta(t)\dot{x}(t),\lambda(t)),x(t)
			\right\rangle
			+\frac{2\alpha}{3}t\beta(t)
			\langle\lambda(t)-\kappa,A\dot{x}(t)\rangle.
		\end{aligned}
	\end{equation}
	Applying the strong convexity inequality for
	$f(\cdot)+\frac{\epsilon(t)}{2}\Vert\cdot\Vert^2$
	at $x(t)+\beta(t)\dot{x}(t)$ and $x^*$,
	and using $Ax^*=b$ and the definition of
	$\mathcal{L}_{\rho}$, we obtain, after expanding
	the quadratic terms,
	\begin{equation}\label{DS_lem1_3}
		\begin{aligned}
			&-\left\langle
			\nabla_x\mathcal{L}_{\rho}
			(x(t)+\beta(t)\dot{x}(t),\kappa)
			+\epsilon(t)x(t),x(t)-x^*
			\right\rangle\\
			&\quad\leq
			\mathcal{L}_{\rho}(x^*,\kappa)
			-\mathcal{L}_{\rho}(x(t)+\beta(t)\dot{x}(t),\kappa)
			-\frac{\rho}{2}
			\Vert A(x(t)+\beta(t)\dot{x}(t))-b\Vert^2\\
			&\qquad
			+\beta(t)\langle\kappa-\lambda(t),A\dot{x}(t)\rangle
			+\frac{\epsilon(t)}{2}
			\left(
			\Vert x^*\Vert^2-\Vert x(t)\Vert^2
			-\Vert x(t)-x^*\Vert^2
			\right)\\
			&\qquad
			+\beta(t)\left\langle
			\nabla_x\mathcal{L}_{\rho}
			(x(t)+\beta(t)\dot{x}(t),\lambda(t)),\dot{x}(t)
			\right\rangle.
		\end{aligned}
	\end{equation}
	Moreover, Young's inequality gives
	\begin{equation}\label{DS_lem1_4}
		\begin{aligned}
			&t\left(t\dot{\beta}(t)-\frac{\alpha}{3}\beta(t)\right)
			\left\langle
			\nabla_x\mathcal{L}_{\rho}
			(x(t)+\beta(t)\dot{x}(t),\lambda(t)),\dot{x}(t)
			\right\rangle\\
			&\quad\leq
			\frac{1}{2}
			\left\vert
			\dot{\beta}(t)-\frac{\alpha\beta(t)}{3t}
			\right\vert
			\left(
			t^{\frac{5}{2}}
			\left\Vert\nabla_x\mathcal{L}_{\rho}
			(x(t)+\beta(t)\dot{x}(t),\lambda(t))\right\Vert^2
			+t^{\frac{3}{2}}\Vert\dot{x}(t)\Vert^2
			\right),
		\end{aligned}
	\end{equation}
	and
	\begin{equation}\label{DS_lem1_5}
		\begin{aligned}
			&-t^2\beta(t)\epsilon(t)
			\left\langle
			\nabla_x\mathcal{L}_{\rho}
			(x(t)+\beta(t)\dot{x}(t),\lambda(t)),x(t)
			\right\rangle\\
			&\quad\leq
			\frac{|\beta(t)|t^2}{M}
			\left\Vert\nabla_x\mathcal{L}_{\rho}
			(x(t)+\beta(t)\dot{x}(t),\lambda(t))\right\Vert^2
			+\frac{M|\beta(t)|t^2\epsilon(t)^2}{4}
			\Vert x(t)\Vert^2,
		\end{aligned}
	\end{equation}
	where $M>1$ is the constant in~\eqref{parameter_condition_1}.
	Substituting~\eqref{DS_lem1_3}--\eqref{DS_lem1_5}
	into~\eqref{DS_lem1_2}, cancelling the terms involving
	$\langle\lambda(t)-\kappa,A\dot{x}(t)\rangle$,
	and discarding the nonpositive penalty term yield
	\begin{equation}\label{DS_lem1_7}
		\begin{aligned}
			\dot{\mathcal{E}}_{\kappa}(t)
			\leq{}&
			\frac{2(3-\alpha)}{3}t
			\left(
			\mathcal{L}_{\rho}(x(t)+\beta(t)\dot{x}(t),\kappa)
			-\mathcal{L}_{\rho}(x^*,\kappa)
			\right)\\
			&-\left(
			\left(\beta(t)-\frac{|\beta(t)|}{M}\right)t^2
			-\frac{t^{\frac{5}{2}}}{2}
			\left\vert
			\dot{\beta}(t)-\frac{\alpha\beta(t)}{3t}
			\right\vert
			\right)
			\left\Vert\nabla_x\mathcal{L}_{\rho}
			(x(t)+\beta(t)\dot{x}(t),\lambda(t))\right\Vert^2\\
			&+\left(
			\frac{3-\alpha}{3}t\epsilon(t)
			+\frac{t^2\dot{\epsilon}(t)}{2}
			+\frac{M|\beta(t)|t^2\epsilon(t)^2}{4}
			\right)\Vert x(t)\Vert^2\\
			&+\frac{\alpha}{3}t\epsilon(t)
			\left(\Vert x^*\Vert^2-\Vert x(t)-x^*\Vert^2\right)\\
			&+\left(
			\frac{3-\alpha}{3}t
			+\frac{t^{\frac{3}{2}}}{2}
			\left\vert
			\dot{\beta}(t)-\frac{\alpha\beta(t)}{3t}
			\right\vert
			\right)\Vert\dot{x}(t)\Vert^2.
		\end{aligned}
	\end{equation}
	Since $M>1$, $\alpha>3$, and
	$\beta(t)=\gamma+\frac{\beta}{t}$ with either
	$\gamma>0$ and $\beta\in\mathbb{R}$, or
	$\gamma=0$ and $\beta>0$,
	there exist $t_1'\geq t_0$ and $r_1>0$ such that,
	for all $t\geq t_1'$,
	\[
	-\left(
	\left(\beta(t)-\frac{|\beta(t)|}{M}\right)t^2
	-\frac{t^{\frac{5}{2}}}{2}
	\left\vert
	\dot{\beta}(t)-\frac{\alpha\beta(t)}{3t}
	\right\vert
	\right)\leq0
	\]
	and
	\[
	\frac{3-\alpha}{3}t
	+\frac{t^{\frac{3}{2}}}{2}
	\left\vert
	\dot{\beta}(t)-\frac{\alpha\beta(t)}{3t}
	\right\vert
	\leq-r_1t.
	\]
	Furthermore, condition~\eqref{parameter_condition_1} implies
	\[
	\frac{3-\alpha}{3}t\epsilon(t)
	+\frac{t^2\dot{\epsilon}(t)}{2}
	+\frac{M|\beta(t)|t^2\epsilon(t)^2}{4}
	\leq\frac{3-\alpha}{3}t\epsilon(t),
	\qquad t\geq t_1.
	\]
	Substituting these estimates into~\eqref{DS_lem1_7}
	yields~\eqref{DS_lem1_6} for all
	$t\geq t_2:=\max\{t_1,t_1'\}$.
\end{proof}

We next recall a basic and widely used technical lemma, which will be instrumental in the subsequent convergence rate analysis of the objective residual and the feasibility violation.
\begin{lemma}{\rm{\cite[Lemma A.3]{Attouch JMAA 2018}}}\label{lemB.2}
	Suppose that $\delta>0$ and $\phi\in L^1\left(\delta,+\infty\right)$ is a non-negative and continuous function, and $\psi:\left[\delta,+\infty\right)\rightarrow\left(0,+\infty\right)$ is a non-decreasing function such that $\lim\limits_{t\rightarrow+\infty}\psi(t)=+\infty$. Then,
	\begin{equation*}
		\lim\limits_{t\rightarrow+\infty}\frac{1}{\psi(t)}\int_{\delta}^{t}\psi(s)\phi(s)ds =0.
	\end{equation*}
\end{lemma}

Based on the energy estimate established in the above lemmas, we proceed to analyze the asymptotic convergence properties of the dynamical system \eqref{DS}.

\begin{theorem}\label{DS_thm1}
	Suppose that the assumptions of Lemma~\ref{DS_lem1} hold
	and that
	$\int_{t_0}^{+\infty}\frac{\epsilon(t)}{t}\,dt<+\infty$.
	Let $(x(t),\lambda(t))$ be a trajectory of the dynamical
	system~\eqref{DS}, and let $(x^*,\lambda^*)\in\Omega$.
	Then $\dot{x}(t)\to\mathbf{0}$ as $t\to+\infty$, and
	\[
	\begin{aligned}
		\mathcal{L}_{\rho}
		\left(x(t)+\left(\gamma+\frac{\beta}{t}\right)\dot{x}(t),
		\lambda^*\right)
		-\mathcal{L}_{\rho}(x^*,\lambda^*)\to0,\quad
	    A\left(x(t)+\left(\gamma+\frac{\beta}{t}\right)
		\dot{x}(t)\right)-b\to\mathbf{0}.
	\end{aligned}
	\]
	If, in addition,
	$\int_{t_0}^{+\infty}t\epsilon(t)\,dt<+\infty$,
	then $x(\cdot)$ is bounded and the following estimates hold:
	\begin{itemize}
		\item[(a)]
		$\mathcal{L}_{\rho}
		\left(x(t)+\left(\gamma+\frac{\beta}{t}\right)\dot{x}(t),
		\lambda^*\right)
		-\mathcal{L}_{\rho}(x^*,\lambda^*)
		=\mathcal{O}\left(\frac{1}{t^2}\right)$,
		\quad
		$\Vert\dot{x}(t)\Vert
		=\mathcal{O}\left(\frac{1}{t}\right)$;
		
		\item[(b)]
		$\left\Vert
		A\left(x(t)+\left(\gamma+\frac{\beta}{t}\right)
		\dot{x}(t)\right)-b
		\right\Vert
		=\mathcal{O}\left(\frac{1}{t^2}\right)$,
		\quad
		$\left|
		f\left(x(t)+\left(\gamma+\frac{\beta}{t}\right)
		\dot{x}(t)\right)-f(x^*)
		\right|
		=\mathcal{O}\left(\frac{1}{t^2}\right)$.
	\end{itemize}
\end{theorem}

\begin{proof}
	Since $(x^*,\lambda^*)\in\Omega$, convexity of $f$ and
	the KKT conditions imply
	\[
	\mathcal{L}_{\rho}(x(t)+\beta(t)\dot{x}(t),\lambda^*)
	-\mathcal{L}_{\rho}(x^*,\lambda^*)
	\geq
	\frac{\rho}{2}
	\Vert A(x(t)+\beta(t)\dot{x}(t))-b\Vert^2
	\geq0.
	\]
	Taking $\kappa=\lambda^*$ in~\eqref{DS_lem1_6},
	discarding the nonpositive terms, and integrating over
	$[t_2,t]$ give
	\begin{equation}\label{DS_thm1_ex1}
		\mathcal{E}_{\lambda^*}(t)
		\leq
		\mathcal{E}_{\lambda^*}(t_2)
		+\frac{\alpha}{3}\Vert x^*\Vert^2
		\int_{t_2}^{t}s\epsilon(s)\,ds.
	\end{equation}
	Applying Lemma~\ref{lemB.2} with
	$\psi(t)=t^2$ and $\phi(t)=\epsilon(t)/t$ yields
	\[
	\frac{1}{t^2}\int_{t_2}^{t}s\epsilon(s)\,ds\to0,
	\qquad
	\frac{\mathcal{E}_{\lambda^*}(t)}{t^2}\to0.
	\]
	By the definition of $\mathcal{E}_{\lambda^*}$ and
	the nonnegativity of its terms, we obtain
	\[
	\begin{aligned}
		&\mathcal{L}_{\rho}(x(t)+\beta(t)\dot{x}(t),\lambda^*)
		-\mathcal{L}_{\rho}(x^*,\lambda^*)\to0,\\
		&\left\Vert
		\frac{2\alpha}{3t}(x(t)-x^*)+\dot{x}(t)
		\right\Vert\to0,
		\qquad
		\frac{\Vert x(t)-x^*\Vert}{t}\to0.
	\end{aligned}
	\]
	Thus $\dot{x}(t)\to\mathbf{0}$, and the first inequality
	of the proof gives
	$A(x(t)+\beta(t)\dot{x}(t))-b\to\mathbf{0}$.
	
	Now assume additionally that
	$\int_{t_0}^{+\infty}t\epsilon(t)\,dt<+\infty$.
	Estimate~\eqref{DS_thm1_ex1} shows that
	$\mathcal{E}_{\lambda^*}$ is bounded on $[t_0,+\infty)$.
	Its definition therefore yields the boundedness of
	$x(\cdot)$ and $\lambda(\cdot)$, together with
	\[
	\begin{aligned}
		&\mathcal{L}_{\rho}(x(t)+\beta(t)\dot{x}(t),\lambda^*)
		-\mathcal{L}_{\rho}(x^*,\lambda^*)
		=\mathcal{O}\left(\frac{1}{t^2}\right),\\
		&\sup_{t\geq t_0}
		\left\Vert
		\frac{2\alpha}{3}(x(t)-x^*)+t\dot{x}(t)
		\right\Vert<+\infty.
	\end{aligned}
	\]
	Consequently,
	$\Vert\dot{x}(t)\Vert=\mathcal{O}(1/t)$,
	which proves~{\rm(a)}.
	
	To prove~{\rm(b)}, for each fixed $s\geq t_2$, define
	\[
	\kappa(s):=
	\begin{cases}
		\displaystyle
		\lambda^*
		+\frac{A(x(s)+\beta(s)\dot{x}(s))-b}
		{\Vert A(x(s)+\beta(s)\dot{x}(s))-b\Vert},
		&\text{if }A(x(s)+\beta(s)\dot{x}(s))-b
		\neq\mathbf{0},\\[1ex]
		\lambda^*,&\text{otherwise}.
	\end{cases}
	\]
	Then $\Vert\kappa(s)-\lambda^*\Vert\leq1$.
	Using~\eqref{DS_lem1_6}, the definition of
	$\mathcal{L}_{\rho}$, and the nonnegativity of the
	augmented Lagrangian gap at $\lambda^*$, we obtain
	\[
	\begin{aligned}
		\dot{\mathcal{E}}_{\kappa(s)}(t)
		\leq{}&
		\frac{2(3-\alpha)}{3}t
		\langle\kappa(s)-\lambda^*,
		A(x(t)+\beta(t)\dot{x}(t))-b\rangle
		+\frac{\alpha}{3}t\epsilon(t)\Vert x^*\Vert^2.
	\end{aligned}
	\]
	Multiplying by $t^{\frac{2(\alpha-3)}{3}}$ and
	integrating by parts give
	\begin{equation}\label{DS_thm1_1}
		\begin{aligned}
			&t^{\frac{2(\alpha-3)}{3}}\mathcal{E}_{\kappa(s)}(t)
			-t_2^{\frac{2(\alpha-3)}{3}}\mathcal{E}_{\kappa(s)}(t_2)
			-\frac{2(\alpha-3)}{3}
			\int_{t_2}^{t}
			\tau^{\frac{2(\alpha-3)}{3}-1}
			\mathcal{E}_{\kappa(s)}(\tau)\,d\tau\\
			&\quad\leq
			\frac{2(3-\alpha)}{3}
			\int_{t_2}^{t}
			\tau^{\frac{2(\alpha-3)}{3}+1}
			\langle\kappa(s)-\lambda^*,
			A(x(\tau)+\beta(\tau)\dot{x}(\tau))-b\rangle\,d\tau\\
			&\qquad
			+\frac{\alpha}{3}\Vert x^*\Vert^2
			\int_{t_2}^{t}
			\tau^{\frac{2(\alpha-3)}{3}+1}
			\epsilon(\tau)\,d\tau.
		\end{aligned}
	\end{equation}
	By $\Vert\kappa(s)-\lambda^*\Vert\leq1$ and
	the definition of $\mathcal{E}_{\lambda^*}$,
	\begin{equation}\label{DS_thm1_2}
		\begin{aligned}
			&\frac{1}{2}
			\left\Vert
			\frac{2\alpha}{3}(x(t)-x^*)+t\dot{x}(t)
			\right\Vert^2
			+\frac{\alpha(\alpha-3)}{9}\Vert x(t)-x^*\Vert^2
			+\frac{\alpha}{3}\Vert\lambda(t)-\kappa(s)\Vert^2\\
			&\quad\leq
			\mathcal{E}_{\lambda^*}(t)
			+\frac{\alpha}{3}\Vert\lambda(t)-\lambda^*\Vert^2
			+\frac{2\alpha}{3}.
		\end{aligned}
	\end{equation}
	It follows that
	\[
	\begin{aligned}
		\mathcal{E}_{\kappa(s)}(t)
		\leq{}&
		2\mathcal{E}_{\lambda^*}(t)+\frac{2\alpha}{3}
		+t^2\langle\kappa(s)-\lambda^*,
		A(x(t)+\beta(t)\dot{x}(t))-b\rangle\\
		\leq{}&
		C_1+t^2\langle\kappa(s)-\lambda^*,
		A(x(t)+\beta(t)\dot{x}(t))-b\rangle,
	\end{aligned}
	\]
	where $C_1>0$ is independent of $s$ and $t\geq t_2$.
	Substituting this bound into~\eqref{DS_thm1_1}
	cancels the integrals involving the constraint residual
	and yields
	\[
	\begin{aligned}
		t^{\frac{2(\alpha-3)}{3}}\mathcal{E}_{\kappa(s)}(t)
		\leq{}&
		t_2^{\frac{2(\alpha-3)}{3}}\mathcal{E}_{\kappa(s)}(t_2)
		+C_1\left(
		t^{\frac{2(\alpha-3)}{3}}
		-t_2^{\frac{2(\alpha-3)}{3}}
		\right)\\
		&+\frac{\alpha}{3}\Vert x^*\Vert^2
		\int_{t_2}^{t}
		\tau^{\frac{2(\alpha-3)}{3}+1}
		\epsilon(\tau)\,d\tau.
	\end{aligned}
	\]
	Hence
	\[
	\begin{aligned}
		\mathcal{E}_{\kappa(s)}(t)
		\leq{}&
		\left(\frac{t_2}{t}\right)^{\frac{2(\alpha-3)}{3}}
		\mathcal{E}_{\kappa(s)}(t_2)+C_1\\
		&+\frac{\alpha}{3}\Vert x^*\Vert^2
		t^{\frac{2(3-\alpha)}{3}}
		\int_{t_2}^{t}
		\tau^{\frac{2(\alpha-3)}{3}+1}
		\epsilon(\tau)\,d\tau.
	\end{aligned}
	\]
	Lemma~\ref{lemB.2}, applied with
	$\psi(t)=t^{\frac{2(\alpha-3)}{3}}$ and
	$\phi(t)=t\epsilon(t)$, shows that the last term
	tends to zero and is therefore bounded by a constant
	$C_2>0$ for $t\geq t_2$.
	Moreover, evaluating the preceding bound for
	$\mathcal{E}_{\kappa(s)}$ at $t_2$ gives
	\[
	\mathcal{E}_{\kappa(s)}(t_2)
	\leq
	C_1+t_2^2
	\Vert A(x(t_2)+\beta(t_2)\dot{x}(t_2))-b\Vert
	=:C_3.
	\]
	Thus
	\[
	\mathcal{E}_{\kappa(s)}(t)
	\leq C_5:=C_1+C_2+C_3,
	\qquad t,s\geq t_2.
	\]
	In particular, the bound is uniform with respect to $s$.
	
	By the definition of $\mathcal{E}_{\kappa(s)}$,
	discarding its nonnegative quadratic terms gives
	\[
	t^2\left(
	f(x(t)+\beta(t)\dot{x}(t))-f(x^*)
	+\langle\kappa(s),
	A(x(t)+\beta(t)\dot{x}(t))-b\rangle
	\right)\leq C_5.
	\]
	Setting $t=s$ and using the definition of $\kappa(s)$,
	we obtain
	\begin{equation}\label{DS_thm1_3}
		\begin{aligned}
			s^2\bigl(
			&\mathcal{L}(x(s)+\beta(s)\dot{x}(s),\lambda^*)
			-\mathcal{L}(x^*,\lambda^*)\\
			&+\Vert A(x(s)+\beta(s)\dot{x}(s))-b\Vert
			\bigr)\leq C_5,
		\end{aligned}
	\end{equation}
	where
	$\mathcal{L}(x,\lambda):=f(x)+\langle\lambda,Ax-b\rangle$.
	The KKT conditions imply
	$\mathcal{L}(x(s)+\beta(s)\dot{x}(s),\lambda^*)
	-\mathcal{L}(x^*,\lambda^*)\geq0$.
	Consequently,
	\[
	\Vert A(x(t)+\beta(t)\dot{x}(t))-b\Vert
	\leq\frac{C_5}{t^2},
	\qquad t\geq t_2.
	\]
	Finally, \eqref{DS_thm1_3} and the Cauchy--Schwarz
	inequality give
	\[
	\begin{aligned}
		f(x(t)+\beta(t)\dot{x}(t))-f(x^*)
		&\leq
		\frac{C_5}{t^2}
		-\langle\lambda^*,
		A(x(t)+\beta(t)\dot{x}(t))-b\rangle\\
		&\leq
		\frac{(1+\Vert\lambda^*\Vert)C_5}{t^2}.
	\end{aligned}
	\]
	On the other hand, convexity of $f$ and
	$(x^*,\lambda^*)\in\Omega$ yield
	\[
	\begin{aligned}
		f(x(t)+\beta(t)\dot{x}(t))-f(x^*)
		&\geq
		-\langle\lambda^*,
		A(x(t)+\beta(t)\dot{x}(t))-b\rangle\\
		&\geq-\frac{\Vert\lambda^*\Vert C_5}{t^2}.
	\end{aligned}
	\]
	Together, these inequalities prove
	$\left|f(x(t)+\beta(t)\dot{x}(t))-f(x^*)\right|
	=\mathcal{O}(1/t^2)$, completing the proof.
\end{proof}
\begin{remark}
	He et al.~\cite{He Optim 2025} proposed a mixed primal--dual dynamical system with an explicit Hessian-driven damping term and established an $\mathcal{O}(1/t^{2})$ convergence rate. This rate is consistent with the $\mathcal{O}(1/t^{2})$ estimate obtained in Theorem \ref{DS_thm1} for our proposed dynamics \eqref{DS}. A key difference is that our dynamical system \eqref{DS} attains this accelerated rate without assuming twice differentiability of the objective function. More recently, Li et al.~\cite{Li arXiv 2025} introduced a fully second-order primal--dual dynamical system and proved that it also enjoys an $\mathcal{O}(1/t^{2})$ rate. In contrast, in our model the dual variable evolves according to a first-order ordinary differential equation. From a numerical perspective, this mixed ``second-order''+``first-order'' structure can substantially reduce the computational cost.
\end{remark}
\begin{remark}\label{rem:critical_epsilon}
	The choice $\epsilon(t)=a/t^r$ with $a>0$ and $r>2$ satisfies the assumptions of Theorem~\ref{DS_thm1}. In contrast, the critical case $\epsilon(t)=a/t^2$ does not satisfy the condition $\int_{t_0}^{+\infty}t\epsilon(t)\,dt<+\infty$. Hence, Theorem~\ref{DS_thm1} does not directly yield the fast convergence rates for this critical case. We next show that these rates can still be obtained for $\epsilon(t)=a/t^2$ by a separate analysis.
\end{remark}

We next establish convergence of the dual trajectory, which
will be used in the critical regularization case and the
minimum-norm convergence analysis.

\begin{lemma}
	\label{DS_lemma3}
	Suppose that the assumptions of Lemma~\ref{DS_lem1} hold
	and that
	$\int_{t_0}^{+\infty}\frac{\epsilon(t)}{t}\,dt<+\infty$.
	Let $(x(t),\lambda(t))$ be a trajectory of the dynamical
	system~\eqref{DS}. Then there exists
	$(x^*,\lambda^\infty)\in\Omega$ such that
	\[
	\lambda(t)\to\lambda^\infty
	\qquad\text{as }t\to+\infty.
	\]
\end{lemma}

\begin{proof}
	If $\operatorname{Ran}A=\{0\}$, feasibility gives $b=0$,
	and the dual equation in~\eqref{DS} yields
	$\dot{\lambda}(t)=0$.
	Since $A^\top=0$, the constant multiplier forms a
	primal--dual solution with any primal minimizer.
	Hence we assume that $\operatorname{Ran}A\neq\{0\}$.
	
	Take any $(x^*,\bar{\lambda})\in\Omega$ and recall that
	$\hat{x}(t)=x(t)+\beta(t)\dot{x}(t)$.
	By Theorem~\ref{DS_thm1},
	\[
	\begin{aligned}
		\dot{x}(t)\to\mathbf{0},\qquad
		A\hat{x}(t)-b\to\mathbf{0},\qquad
		\mathcal{L}_{\rho}(\hat{x}(t),\bar{\lambda})
		-\mathcal{L}_{\rho}(x^*,\bar{\lambda})\to0.
	\end{aligned}
	\]
	Since $\beta(t)$ is bounded, we also have
	$\hat{x}(t)-x(t)\to\mathbf{0}$ and $Ax(t)-b\to\mathbf{0}$.
	Moreover, convexity of $f$, Lipschitz continuity of
	$\nabla f$, and the KKT conditions imply
	\[
	\begin{aligned}
		\Vert\nabla f(\hat{x}(t))-\nabla f(x^*)\Vert^2
		&\leq
		2L\bigl(
		f(\hat{x}(t))-f(x^*)
		-\langle\nabla f(x^*),\hat{x}(t)-x^*\rangle
		\bigr)\\
		&=
		2L\left(
		\mathcal{L}_{\rho}(\hat{x}(t),\bar{\lambda})
		-\mathcal{L}_{\rho}(x^*,\bar{\lambda})
		-\frac{\rho}{2}\Vert A\hat{x}(t)-b\Vert^2
		\right)\to0.
	\end{aligned}
	\]
	
	Every term in the dual equation belongs to
	$\operatorname{Ran}A$, so
	$P_{\ker A^\top}\dot{\lambda}(t)=0$.
	Define
	\[
	\lambda^\infty
	:=P_{\operatorname{Ran}A}\bar{\lambda}
	+P_{\ker A^\top}\lambda(t_0).
	\]
	Then
	$A^\top\lambda^\infty=A^\top\bar{\lambda}$,
	hence $(x^*,\lambda^\infty)\in\Omega$, and
	$\lambda(t)-\lambda^\infty\in\operatorname{Ran}A$
	for all $t\geq t_0$.
	We now distinguish two cases.
	
	\noindent
	\textbf{Case I: $\gamma>0$ and $\beta\in\mathbb{R}$.}
	
	There exists $T\geq t_0$ such that
	$\beta(t)\geq\gamma/2>0$ for $t\geq T$.
	Set $e(t):=\lambda(t)-\lambda^\infty$.
	The dual equation in~\eqref{DS} gives
	\[
	\dot{e}(t)+c(t)AA^\top e(t)=q(t),
	\qquad
	c(t):=\frac{3t^2\beta(t)}{2\alpha},
	\]
	where
	\[
	\begin{aligned}
		\frac{q(t)}{c(t)}
		={}&
		\frac{2\alpha}{3t\beta(t)}(Ax(t)-b)
		-\epsilon(t)Ax(t)
		+\frac{\eta(t)}{3t\beta(t)}A\dot{x}(t)\\
		&-A\bigl(\nabla f(\hat{x}(t))-\nabla f(x^*)\bigr)
		-\rho AA^\top(A\hat{x}(t)-b).
	\end{aligned}
	\]
	Since $\beta(t)\to\gamma>0$, the coefficients
	$1/(t\beta(t))$ and $\eta(t)/(3t\beta(t))$ are bounded
	on $[T,+\infty)$.
	Also, $Ax(t)\to b$ and $\epsilon(t)\to0$ imply
	$\epsilon(t)Ax(t)\to\mathbf{0}$.
	The preceding convergence results therefore give
	$\Vert q(t)\Vert/c(t)\to0$.
	
	Since $AA^\top$ is positive definite on
	$\operatorname{Ran}A$, there exists $\tilde{\sigma}>0$
	such that
	\[
	\langle AA^\top y,y\rangle
	\geq\tilde{\sigma}\Vert y\Vert^2,
	\qquad y\in\operatorname{Ran}A.
	\]
	Taking the inner product of the equation for $e(t)$
	with $e(t)$ and applying Young's inequality yield,
	for $w(t):=\Vert e(t)\Vert^2$,
	\[
	w'(t)+\tilde{\sigma}c(t)w(t)
	\leq
	\frac{1}{\tilde{\sigma}}
	\frac{\Vert q(t)\Vert^2}{c(t)}.
	\]
	For any $\tilde{\delta}>0$, choose
	$T_{\tilde{\delta}}\geq T$ such that
	$\Vert q(t)\Vert/c(t)\leq\tilde{\delta}$
	for $t\geq T_{\tilde{\delta}}$.
	Then
	\[
	w'(t)+\tilde{\sigma}c(t)
	\left(w(t)-\frac{\tilde{\delta}^2}{\tilde{\sigma}^2}\right)
	\leq0.
	\]
	Multiplying by the integrating factor
	$\exp\bigl(\tilde{\sigma}
	\int_{T_{\tilde{\delta}}}^{t}c(s)\,ds\bigr)$
	and integrating give
	\[
	\begin{aligned}
		w(t)\leq{}&
		\frac{\tilde{\delta}^2}{\tilde{\sigma}^2}
		+\left[
		w(T_{\tilde{\delta}})
		-\frac{\tilde{\delta}^2}{\tilde{\sigma}^2}
		\right]_+
		\exp\left(
		-\tilde{\sigma}
		\int_{T_{\tilde{\delta}}}^{t}c(s)\,ds
		\right),
	\end{aligned}
	\]
	where $[r]_+:=\max\{r,0\}$.
	Since $c(t)\sim\frac{3\gamma}{2\alpha}t^2$,
	the integral in the exponent diverges, and thus
	\[
	\limsup_{t\to+\infty}w(t)
	\leq\frac{\tilde{\delta}^2}{\tilde{\sigma}^2}.
	\]
	Letting $\tilde{\delta}\downarrow0$ proves
	$e(t)\to\mathbf{0}$ and hence
	$\lambda(t)\to\lambda^\infty$.
	
	\noindent
	\textbf{Case II: $\gamma=0$ and $\beta>0$.}
	
	Here $\beta(t)=\beta/t$.
	Set $u(t):=\lambda(t)-\lambda^\infty$ and let
	\[
	\mathsf{H}:=
	\left.AA^\top\right|_{\operatorname{Ran}A},
	\]
	which is positive definite on $\operatorname{Ran}A$.
	Define
	\[
	\begin{aligned}
		\zeta(t):={}&
		A\bigl(\nabla f(\hat{x}(t))-\nabla f(x^*)\bigr)
		+\rho\mathsf{H}(A\hat{x}(t)-b)
		+\epsilon(t)Ax(t).
	\end{aligned}
	\]
	The convergence results established above imply
	$\zeta(t)\to\mathbf{0}$.
	Using the primal equation and
	$\nabla f(x^*)+A^\top\lambda^\infty=0$, we obtain
	\begin{equation}\label{ex_DS_lem2.3_1}
		A\ddot{x}(t)+\frac{\alpha}{t}A\dot{x}(t)
		=-\mathsf{H}u(t)-\zeta(t).
	\end{equation}
	Define
	\[
	P(t):=\frac{2\alpha}{3}(Ax(t)-b)+tA\dot{x}(t).
	\]
	Differentiating and using~\eqref{ex_DS_lem2.3_1} give
	\begin{equation}\label{ex_DS_lem2.3_2}
		\begin{aligned}
			\dot{P}(t)
			&=
			\frac{3-\alpha}{3}A\dot{x}(t)
			-t\mathsf{H}u(t)-t\zeta(t)\\
			&=
			-t\mathsf{H}u(t)-t\zeta(t)
			-\frac{\alpha-3}{3t}P(t)
			+\frac{2\alpha(\alpha-3)}{9t}(Ax(t)-b).
		\end{aligned}
	\end{equation}
	Moreover,
	$\eta(t)=3t-(\alpha+3)\beta/t$ and
	\[
	A\dot{x}(t)
	=\frac{1}{t}P(t)-\frac{2\alpha}{3t}(Ax(t)-b).
	\]
	Substituting these identities into the dual equation
	in~\eqref{DS} yields
	\begin{equation}\label{ex_DS_lem2.3_3}
		\begin{aligned}
			\dot{u}(t)
			={}&
			\frac{3t}{2\alpha}
			\bigl(P(t)-\beta\mathsf{H}u(t)-\beta\zeta(t)\bigr)
			-\frac{\beta(\alpha+3)}{2\alpha t}P(t)\\
			&+\frac{\beta(\alpha+3)}{3t}(Ax(t)-b).
		\end{aligned}
	\end{equation}
	For
	\[
	Z(t):=
	\begin{pmatrix}
		P(t)\\
		u(t)
	\end{pmatrix}
	\in\operatorname{Ran}A\times\operatorname{Ran}A,
	\]
	equations~\eqref{ex_DS_lem2.3_2}
	and~\eqref{ex_DS_lem2.3_3} become
	\[
	\dot{Z}(t)
	=t\mathcal{M}Z(t)+t\mathcal{N}\zeta(t)
	+\mathcal{R}_0(t),
	\]
	where
	\[
	\mathcal{M}=
	\begin{pmatrix}
		0&-\mathsf{H}\\
		\frac{3}{2\alpha}I&-\frac{3\beta}{2\alpha}\mathsf{H}
	\end{pmatrix},
	\qquad
	\mathcal{N}=
	\begin{pmatrix}
		-I\\
		-\frac{3\beta}{2\alpha}I
	\end{pmatrix},
	\]
	and
	\[
	\mathcal{R}_0(t)=
	\begin{pmatrix}
		\dfrac{3-\alpha}{3t}P(t)
		+\dfrac{2\alpha(\alpha-3)}{9t}(Ax(t)-b)\\[2mm]
		-\dfrac{\beta(\alpha+3)}{2\alpha t}P(t)
		+\dfrac{\beta(\alpha+3)}{3t}(Ax(t)-b)
	\end{pmatrix}.
	\]
	In particular, there exists $C_0>0$ such that
	\[
	\Vert\mathcal{R}_0(t)\Vert
	\leq\frac{C_0}{t}
	\bigl(\Vert Z(t)\Vert+\Vert Ax(t)-b\Vert\bigr).
	\]
	
	For each eigenvalue $\mu>0$ of $\mathsf{H}$,
	the corresponding block of $\mathcal{M}$ has
	characteristic polynomial
	\[
	\chi_\mu(\xi)
	=\xi^2+\frac{3\beta\mu}{2\alpha}\xi
	+\frac{3\mu}{2\alpha}.
	\]
	Both roots have negative real parts since
	$\alpha,\beta,\mu>0$.
	Thus $\mathcal{M}$ is Hurwitz, and there exists
	a symmetric positive definite operator $\mathcal{S}$
	satisfying
	\[
	\mathcal{M}^\top\mathcal{S}
	+\mathcal{S}\mathcal{M}=-I.
	\]
	Set $V(t):=\langle\mathcal{S}Z(t),Z(t)\rangle$.
	There exist $m_{\mathcal{S}},M_{\mathcal{S}}>0$ such that
	\[
	m_{\mathcal{S}}\Vert Z(t)\Vert^2
	\leq V(t)\leq
	M_{\mathcal{S}}\Vert Z(t)\Vert^2.
	\]
	Differentiating $V(t)$ gives
	\[
	\dot{V}(t)
	=-t\Vert Z(t)\Vert^2
	+2t\langle\mathcal{S}Z(t),\mathcal{N}\zeta(t)\rangle
	+2\langle\mathcal{S}Z(t),\mathcal{R}_0(t)\rangle.
	\]
	By Young's inequality and the bound for
	$\mathcal{R}_0(t)$, there exist $C_1,C_2>0$
	such that, for all sufficiently large $t$,
	\[
	\begin{aligned}
		2t\left|
		\langle\mathcal{S}Z(t),\mathcal{N}\zeta(t)\rangle
		\right|
		&\leq\frac{t}{4}\Vert Z(t)\Vert^2
		+C_1t\Vert\zeta(t)\Vert^2,\\
		2\left|
		\langle\mathcal{S}Z(t),\mathcal{R}_0(t)\rangle
		\right|
		&\leq\frac{t}{4}\Vert Z(t)\Vert^2
		+\frac{C_2}{t^3}\Vert Ax(t)-b\Vert^2.
	\end{aligned}
	\]
	Consequently, there exist
	$\hat{c}_0,C_3,C_4>0$ and $T_1\geq t_0$ such that
	\[
	\dot{V}(t)+\hat{c}_0tV(t)\leq t\,h(t),
	\qquad t\geq T_1,
	\]
	where
	\[
	h(t):=C_3\Vert\zeta(t)\Vert^2
	+\frac{C_4}{t^4}\Vert Ax(t)-b\Vert^2
	\longrightarrow0.
	\]
	For any $\delta>0$, choose $T_\delta\geq T_1$
	such that $h(t)\leq\delta$ for $t\geq T_\delta$.
	Then
	\[
	\dot{V}(t)+\hat{c}_0t
	\left(V(t)-\frac{\delta}{\hat{c}_0}\right)\leq0.
	\]
	The integrating-factor argument used in Case I,
	with coefficient $\hat{c}_0t$, gives
	\[
	\limsup_{t\to+\infty}V(t)
	\leq\frac{\delta}{\hat{c}_0}.
	\]
	Since $\delta>0$ is arbitrary, $V(t)\to0$.
	Positive definiteness of $\mathcal{S}$ therefore yields
	$P(t)\to\mathbf{0}$ and $u(t)\to\mathbf{0}$.
	Hence $\lambda(t)\to\lambda^\infty$, completing the proof.
\end{proof}
\begin{theorem}
	\label{thm:critical_epsilon}
	Suppose that the assumptions of Lemma~\ref{DS_lem1} hold
	with $\epsilon(t)=a/t^2$, where $a>0$.
	Let $(x(t),\lambda(t))$ be a global solution of the
	dynamical system~\eqref{DS}.
	Then $\lambda(\cdot)$ is bounded on $[t_0,+\infty)$.
	Moreover, for any $(x^*,\lambda^*)\in\Omega$,
	\[
	\begin{aligned}
		&\mathcal{L}_{\rho}
		\left(x(t)+\left(\gamma+\frac{\beta}{t}\right)
		\dot{x}(t),\lambda^*\right)
		-\mathcal{L}_{\rho}(x^*,\lambda^*)
		=\mathcal{O}\left(\frac{1}{t^2}\right),
		\qquad
		\Vert\dot{x}(t)\Vert=\mathcal{O}\left(\frac{1}{t}\right),\\
		&\left\Vert
		A\left(x(t)+\left(\gamma+\frac{\beta}{t}\right)
		\dot{x}(t)\right)-b
		\right\Vert
		=\mathcal{O}\left(\frac{1}{t^2}\right),\quad
		\left|
		f\left(x(t)+\left(\gamma+\frac{\beta}{t}\right)
		\dot{x}(t)\right)-f(x^*)
		\right|
		=\mathcal{O}\left(\frac{1}{t^2}\right).
	\end{aligned}
	\]
\end{theorem}

\begin{proof}
	Since $\beta(t)=\gamma+\beta/t$ is bounded on
	$[t_0,+\infty)$, there exists $C_\beta>0$ such that
	$|\beta(t)|\leq C_\beta$.
	For $\epsilon(t)=a/t^2$, we have
	\[
	2\dot{\epsilon}(t)+M|\beta(t)|\epsilon(t)^2
	\leq-\frac{4a}{t^3}
	+\frac{MC_\beta a^2}{t^4}
	\leq0
	\]
	for all sufficiently large $t$.
	Thus~\eqref{parameter_condition_1} is satisfied.
	Also,
	$\int_{t_0}^{+\infty}\epsilon(t)/t\,dt<+\infty$,
	so Lemma~\ref{DS_lemma3} gives convergence, and hence
	boundedness, of $\lambda(\cdot)$.
	
	Applying
	$\Vert u+v\Vert^2\leq2\Vert u\Vert^2+2\Vert v\Vert^2$
	to~\eqref{DS_lem1_1} yields
	\[
	\begin{aligned}
		\mathcal{E}_{\lambda^*}(t)
		\leq{}&
		t^2\bigl(
		\mathcal{L}_{\rho}(x(t)+\beta(t)\dot{x}(t),\lambda^*)
		-\mathcal{L}_{\rho}(x^*,\lambda^*)
		\bigr)
		+\frac{a}{2}\Vert x(t)\Vert^2\\
		&+\frac{\alpha(5\alpha-3)}{9}\Vert x(t)-x^*\Vert^2
		+t^2\Vert\dot{x}(t)\Vert^2
		+\frac{\alpha}{3}\Vert\lambda(t)-\lambda^*\Vert^2.
	\end{aligned}
	\]
	Choose
	\[
	0<c\leq
	\min\left\{
	\frac{2(\alpha-3)}{3},
	\frac{3a}{5\alpha-3},r_1
	\right\},
	\]
	where $r_1>0$ is given by Lemma~\ref{DS_lem1}.
	Combining the preceding inequality with
	\eqref{DS_lem1_6}, we obtain, for all sufficiently
	large $t$,
	\[
	\begin{aligned}
		\dot{\mathcal{E}}_{\lambda^*}(t)
		+\frac{c}{t}\mathcal{E}_{\lambda^*}(t)
		\leq{}&
		\left(\frac{2(3-\alpha)}{3}+c\right)t
		\bigl(
		\mathcal{L}_{\rho}(x(t)+\beta(t)\dot{x}(t),\lambda^*)
		-\mathcal{L}_{\rho}(x^*,\lambda^*)
		\bigr)\\
		&+\left(
		\frac{ac}{2t}+\frac{(3-\alpha)a}{3t}
		\right)\Vert x(t)\Vert^2+\frac{c\alpha}{3t}\Vert\lambda(t)-\lambda^*\Vert^2
		+\frac{a\alpha}{3t}\Vert x^*\Vert^2\\
		&+\left(
		\frac{\alpha c(5\alpha-3)}{9t}
		-\frac{a\alpha}{3t}
		\right)\Vert x(t)-x^*\Vert^2
		+(c-r_1)t\Vert\dot{x}(t)\Vert^2\\
		\leq{}&
		\frac{c\alpha}{3t}\Vert\lambda(t)-\lambda^*\Vert^2
		+\frac{a\alpha}{3t}\Vert x^*\Vert^2
		\leq\frac{C}{t},
	\end{aligned}
	\]
	where the last inequality follows from the boundedness
	of $\lambda(\cdot)$.
	Multiplying by $t^c$ and integrating over $[T,t]$,
	for a sufficiently large $T\geq t_2$, give
	\[
	\mathcal{E}_{\lambda^*}(t)
	\leq
	\left(\frac{T}{t}\right)^c\mathcal{E}_{\lambda^*}(T)
	+\frac{C}{c}
	\left(1-\left(\frac{T}{t}\right)^c\right).
	\]
	Hence $\mathcal{E}_{\lambda^*}$ is bounded.
	Its definition implies that $x(\cdot)$ is bounded and
	\[
	\begin{aligned}
		&\mathcal{L}_{\rho}(x(t)+\beta(t)\dot{x}(t),\lambda^*)
		-\mathcal{L}_{\rho}(x^*,\lambda^*)
		=\mathcal{O}\left(\frac{1}{t^2}\right),\\
		&\sup_{t\geq T}
		\left\Vert
		\frac{2\alpha}{3}(x(t)-x^*)+t\dot{x}(t)
		\right\Vert<+\infty.
	\end{aligned}
	\]
	Consequently,
	$\Vert\dot{x}(t)\Vert=\mathcal{O}(1/t)$.
	
	It remains to estimate the feasibility violation and
	objective residual.
	For $s\geq T$, let $\kappa(s)$ be defined as in the
	proof of Theorem~\ref{DS_thm1}.
	Since $\mathcal{E}_{\lambda^*}$ is bounded and
	$\Vert\kappa(s)-\lambda^*\Vert\leq1$,
	the same energy comparison gives
	\[
	\mathcal{E}_{\kappa(s)}(t)
	\leq C_1+t^2
	\langle\kappa(s)-\lambda^*,
	A(x(t)+\beta(t)\dot{x}(t))-b\rangle,
	\qquad t,s\geq T,
	\]
	where $C_1>0$ is independent of $t$ and $s$.
	Repeating the weighted integration leading to
	\eqref{DS_thm1_1}, with lower limit $T$, yields
	\[
	\begin{aligned}
		\mathcal{E}_{\kappa(s)}(t)
		\leq{}&
		\left(\frac{T}{t}\right)^{\frac{2(\alpha-3)}{3}}
		\mathcal{E}_{\kappa(s)}(T)+C_1
		+\frac{\alpha}{3}\Vert x^*\Vert^2
		t^{\frac{2(3-\alpha)}{3}}
		\int_T^t
		\tau^{\frac{2(\alpha-3)}{3}+1}
		\epsilon(\tau)\,d\tau.
	\end{aligned}
	\]
	In the present critical case,
	\[
	\begin{aligned}
		t^{\frac{2(3-\alpha)}{3}}
		\int_T^t
		\tau^{\frac{2(\alpha-3)}{3}+1}
		\epsilon(\tau)\,d\tau
		=
		\frac{3a}{2(\alpha-3)}
		\left(
		1-\left(\frac{T}{t}\right)^{\frac{2(\alpha-3)}{3}}
		\right)
		\leq\frac{3a}{2(\alpha-3)}.
	\end{aligned}
	\]
	Moreover,
	\[
	\mathcal{E}_{\kappa(s)}(T)
	\leq
	C_1+T^2
	\Vert A(x(T)+\beta(T)\dot{x}(T))-b\Vert
	=:C_2.
	\]
	Thus
	\[
	\mathcal{E}_{\kappa(s)}(t)
	\leq
	C_1+C_2+\frac{a\alpha}{2(\alpha-3)}
	\Vert x^*\Vert^2
	=:C_3,
	\qquad t,s\geq T.
	\]
	Setting $s=t$ and using the definition of $\kappa(s)$,
	as in~\eqref{DS_thm1_3}, we obtain
	\[
	\begin{aligned}
		t^2\bigl(
		\mathcal{L}(x(t)+\beta(t)\dot{x}(t),\lambda^*)
		-\mathcal{L}(x^*,\lambda^*)
		+\Vert A(x(t)+\beta(t)\dot{x}(t))-b\Vert
		\bigr)\leq C_3.
	\end{aligned}
	\]
	The KKT conditions imply that the Lagrangian gap
	in this inequality is nonnegative.
	Therefore,
	\[
	\Vert A(x(t)+\beta(t)\dot{x}(t))-b\Vert
	\leq\frac{C_3}{t^2}.
	\]
	The same upper and lower objective estimates used
	in Theorem~\ref{DS_thm1} then give
	\[
	\left|f(x(t)+\beta(t)\dot{x}(t))-f(x^*)\right|
	\leq\frac{(1+\Vert\lambda^*\Vert)C_3}{t^2}.
	\]
	This completes the proof.
\end{proof}

In what follows, we further prove the strong convergence of the trajectory generated by the dynamical system~\eqref{DS}. 
We first give an auxiliary lemma, which states that the minimizer of the Tikhonov-regularized augmented Lagrangian converges to the minimum-norm solution as the regularization parameter tends to zero.

To establish minimum-norm strong convergence for the
dynamical system~\eqref{DS}, we first prove an auxiliary
result on Tikhonov regularization of the augmented Lagrangian.

\begin{lemma}\label{lem:tikhonov_AL}
	Let $\rho>0$ and let $(x^*,\lambda^\infty)\in\Omega$,
	where $x^*:=\operatorname{Proj}_{S}0$ and
	$S:=\operatorname*{argmin}\{f(x):Ax=b\}$.
	For each $\epsilon>0$, define
	\[
	x_\epsilon
	:=
	\operatorname*{argmin}_{x\in\mathbb R^n}
	\left\{
	\mathcal L_\rho(x,\lambda^\infty)
	+\frac{\epsilon}{2}\|x\|^2
	\right\}.
	\]
	Then $x_\epsilon$ is uniquely defined and
	$x_\epsilon\to x^*$ as $\epsilon\to0^+$.
\end{lemma}

\begin{proof}
	For each $\epsilon>0$, the function
	$\mathcal L_\rho(\cdot,\lambda^\infty)
	+\frac{\epsilon}{2}\|\cdot\|^2$
	is continuous, coercive, and strongly convex on
	$\mathbb R^n$, so $x_\epsilon$ exists and is unique.
	We next show that
	\[
	\operatorname*{argmin}_{x\in\mathbb R^n}
	\mathcal L_\rho(x,\lambda^\infty)=S.
	\]
	By the convexity of $f$ and
	$(x^*,\lambda^\infty)\in\Omega$, for every
	$x\in\mathbb R^n$,
	\[
	\begin{aligned}
		\mathcal L_\rho(x,\lambda^\infty)
		-\mathcal L_\rho(x^*,\lambda^\infty)
		&=
		f(x)-f(x^*)
		+\langle\lambda^\infty,Ax-b\rangle
		+\frac{\rho}{2}\|Ax-b\|^2\\
		&\geq
		\langle\nabla f(x^*),x-x^*\rangle
		+\langle\lambda^\infty,Ax-b\rangle
		+\frac{\rho}{2}\|Ax-b\|^2\\
		&=
		\frac{\rho}{2}\|Ax-b\|^2
		\geq0.
	\end{aligned}
	\]
	Every $z\in S$ satisfies
	$\mathcal L_\rho(z,\lambda^\infty)
	=\mathcal L_\rho(x^*,\lambda^\infty)=f(x^*)$
	and hence minimizes
	$\mathcal L_\rho(\cdot,\lambda^\infty)$.
	Conversely, if $\bar x$ is a minimizer, the preceding
	inequality and $\rho>0$ imply $A\bar x=b$, and thus
	$f(\bar x)=f(x^*)$.
	This proves the claimed equality of the minimizer sets.
	
	By the definition of $x_\epsilon$,
	\[
	\mathcal L_\rho(x_\epsilon,\lambda^\infty)
	+\frac{\epsilon}{2}\|x_\epsilon\|^2
	\leq
	\mathcal L_\rho(x^*,\lambda^\infty)
	+\frac{\epsilon}{2}\|x^*\|^2.
	\]
	Since $x^*$ minimizes
	$\mathcal L_\rho(\cdot,\lambda^\infty)$, it follows that
	$\|x_\epsilon\|\leq\|x^*\|$ and
	\[
	0\leq
	\mathcal L_\rho(x_\epsilon,\lambda^\infty)
	-\mathcal L_\rho(x^*,\lambda^\infty)
	\leq
	\frac{\epsilon}{2}\|x^*\|^2.
	\]
	Thus $\{x_\epsilon\}$ is bounded, and continuity implies
	that every cluster point $\bar x$ as $\epsilon\to0^+$
	belongs to $S$ and satisfies
	$\|\bar x\|\leq\|x^*\|$.
	Since $x^*$ is the unique minimum-norm element of $S$,
	every such cluster point equals $x^*$.
	Therefore, $x_\epsilon\to x^*$ as $\epsilon\to0^+$.
\end{proof}
\begin{theorem}\label{DS_thm2}
	Suppose that the assumptions of Lemma~\ref{DS_lem1} hold
	and that
	\[
	\int_{t_0}^{+\infty}\frac{\epsilon(t)}{t}\,dt<+\infty,
	\qquad
	\lim_{t\to+\infty}
	t^{\frac{\alpha+3}{3}}\epsilon(t)=+\infty,\quad \text{and}\quad\lim_{t\to+\infty}
	\frac{1}{t^{\frac{\alpha+3}{3}}\epsilon(t)}
	\int_{t_0}^{t}
	s^{\frac{\alpha+3}{3}}\epsilon(s)^2\,ds=0.
	\]
	Assume also that
	$t^2\epsilon(t)\geq\frac{2\alpha(\alpha-3)}{9}$
	for all sufficiently large $t$.
	Let $(x(\cdot),\lambda(\cdot))$ be a global solution
	of the dynamical system~\eqref{DS}, and let
	$x^*=\operatorname{Proj}_S0$.
	Then
	\[
	\liminf_{t\to+\infty}\|x(t)-x^*\|=0.
	\]
	Moreover, if there exists $T\geq t_0$ such that
	$\{x(t):t\geq T\}$ lies entirely in either
	$B(0,\|x^*\|)$ or its complement, then
	\[
	\lim_{t\to+\infty}\|x(t)-x^*\|=0.
	\]
\end{theorem}

\begin{proof}
	Since $\epsilon$ is nonincreasing and
	$\int_{t_0}^{+\infty}\epsilon(t)/t\,dt<+\infty$,
	we have $\epsilon(t)\to0$.
	By Theorem~\ref{DS_thm1} and Lemma~\ref{DS_lemma3},
	$\dot{x}(t)\to0$ and $\lambda(t)\to\lambda^\infty$,
	where $(x^*,\lambda^\infty)\in\Omega$.
	Here we use the fact that any primal solution can be
	paired with a dual solution to form a saddle point.
	In particular,
	\[
	\mathcal{L}_{\rho}
	(x(t)+\beta(t)\dot{x}(t),\lambda^\infty)
	-\mathcal{L}_{\rho}(x^*,\lambda^\infty)\to0.
	\]
	We distinguish the following three cases.
	
	\medskip
	\noindent
	\textbf{Case I:}
	There exists $T\geq t_0$ such that
	$\|x(t)\|\geq\|x^*\|$ for all $t\geq T$.
	
	For $t\geq T$, define
	\[
	\begin{aligned}
		\widetilde{\mathcal{E}}(t)
		={}&t^{\frac{\alpha+3}{3}}
		\Big(
		\mathcal{L}_{\rho}
		(x(t)+\beta(t)\dot{x}(t),\lambda^\infty)
		-\mathcal{L}_{\rho}(x^*,\lambda^\infty)
		+\frac{\epsilon(t)}{2}\|x(t)\|^2
		-\frac{\epsilon(t)}{2}\|x^*\|^2
		\Big)\\
		&+\frac{t^{\frac{\alpha-3}{3}}}{2}
		\left\|
		\frac{2\alpha}{3}(x(t)-x^*)+t\dot{x}(t)
		\right\|^2
		+\frac{\alpha}{3}t^{\frac{\alpha-3}{3}}
		\|\lambda(t)-\lambda^\infty\|^2.
	\end{aligned}
	\]
	This energy is nonnegative in the present case, and
	\[
	\begin{aligned}
		\widetilde{\mathcal{E}}(t)
		={}&t^{\frac{\alpha-3}{3}}
		\mathcal{E}_{\lambda^\infty}(t)
		-t^{\frac{\alpha+3}{3}}
		\frac{\epsilon(t)}{2}\|x^*\|^2
		-\frac{\alpha(\alpha-3)}{9}
		t^{\frac{\alpha-3}{3}}\|x(t)-x^*\|^2.
	\end{aligned}
	\]
	Fix $M>1$ as in~\eqref{parameter_condition_1}.
	Differentiating this identity and applying the convexity
	and Young estimates used in Lemma~\ref{DS_lem1}, with
	an auxiliary constant $q>0$, yield
	\[
	\begin{aligned}
		\dot{\widetilde{\mathcal{E}}}(t)
		\leq{}&
		\left(1-\frac{\alpha}{3}\right)t^{\frac{\alpha}{3}}
		\Big(
		\mathcal{L}_{\rho}
		(x(t)+\beta(t)\dot{x}(t),\lambda^\infty)
		-\mathcal{L}_{\rho}(x^*,\lambda^\infty)
		\Big)\\
		&+\left(
		\frac{|\beta(t)|}{M}-\beta(t)
		+\frac{t}{q}
		\left|\dot{\beta}(t)-\frac{\alpha\beta(t)}{3t}\right|
		\right)t^{\frac{\alpha+3}{3}}
		\left\|
		\nabla_x\mathcal{L}_{\rho}
		(x(t)+\beta(t)\dot{x}(t),\lambda(t))
		\right\|^2\\
		&+\left(
		\frac{3-\alpha}{6}
		+\frac{q}{4}
		\left|\dot{\beta}(t)-\frac{\alpha\beta(t)}{3t}\right|
		\right)t^{\frac{\alpha}{3}}\|\dot{x}(t)\|^2
		-\frac{\rho\alpha}{3}t^{\frac{\alpha}{3}}
		\|A(x(t)+\beta(t)\dot{x}(t))-b\|^2\\
		&+\left(
		-\frac{\alpha}{3}t^{\frac{\alpha}{3}}\epsilon(t)
		+\frac{2(\alpha-3)\alpha^2}{27}
		t^{\frac{\alpha-6}{3}}
		\right)\|x(t)-x^*\|^2\\
		&+t^{\frac{\alpha+3}{3}}
		\left(
		\frac{\dot{\epsilon}(t)}{2}
		+\frac{3-\alpha}{6t}\epsilon(t)
		+\frac{M|\beta(t)|}{4}\epsilon(t)^2
		\right)
		\big(\|x(t)\|^2-\|x^*\|^2\big)\\
		&+\frac{M|\beta(t)|}{4}
		t^{\frac{\alpha+3}{3}}\epsilon(t)^2\|x^*\|^2
		+\frac{\alpha(\alpha-3)}{9}
		t^{\frac{\alpha-6}{3}}
		\|\lambda(t)-\lambda^\infty\|^2.
	\end{aligned}
	\]
	Because $\beta(t)=\gamma+\beta/t$, $M>1$, and
	$\alpha>3$, we can first choose $q>0$ sufficiently
	large and then $t_3\geq T$ such that, for all $t\geq t_3$,
	\[
	\begin{aligned}
		\frac{|\beta(t)|}{M}-\beta(t)
		+\frac{t}{q}
		\left|\dot{\beta}(t)-\frac{\alpha\beta(t)}{3t}\right|
		\leq0\quad\text{ and }\quad
		\frac{3-\alpha}{6}
		+\frac{q}{4}
		\left|\dot{\beta}(t)-\frac{\alpha\beta(t)}{3t}\right|
		\leq0.
	\end{aligned}
	\]
	Condition~\eqref{parameter_condition_1} makes the
	coefficient of $\|x(t)\|^2-\|x^*\|^2$ nonpositive.
	The coefficient of $\|x(t)-x^*\|^2$ is also nonpositive
	by $t^2\epsilon(t)\geq2\alpha(\alpha-3)/9$.
	Discarding the nonpositive terms, we obtain
	\[
	\dot{\widetilde{\mathcal{E}}}(t)
	\leq
	\frac{M|\beta(t)|}{4}
	t^{\frac{\alpha+3}{3}}\epsilon(t)^2\|x^*\|^2
	+\frac{\alpha(\alpha-3)}{9}
	t^{\frac{\alpha-6}{3}}
	\|\lambda(t)-\lambda^\infty\|^2.
	\]
	Since $\beta(t)$ is bounded, integration gives a
	constant $C_6>0$ such that
	\begin{equation}\label{DS_thm2_1}
		\begin{aligned}
			\widetilde{\mathcal{E}}(t)
			\leq{}
			\widetilde{\mathcal{E}}(t_3)
			+C_6\|x^*\|^2
			\int_{t_3}^{t}
			s^{\frac{\alpha+3}{3}}\epsilon(s)^2\,ds
			+\frac{\alpha(\alpha-3)}{9}
			\int_{t_3}^{t}
			s^{\frac{\alpha-6}{3}}
			\|\lambda(s)-\lambda^\infty\|^2\,ds.
		\end{aligned}
	\end{equation}
	
	Let $x_{\epsilon(t)}$ be the regularized minimizer
	defined in Lemma~\ref{lem:tikhonov_AL}.
	Strong convexity and the optimality of $x_{\epsilon(t)}$
	imply, for every $x\in\mathbb R^n$,
	\[
	\begin{aligned}
		\frac{\epsilon(t)}{2}\|x-x_{\epsilon(t)}\|^2
		\leq{}
		\mathcal{L}_{\rho}(x,\lambda^\infty)
		+\frac{\epsilon(t)}{2}\|x\|^2
		-\mathcal{L}_{\rho}(x_{\epsilon(t)},\lambda^\infty)
		-\frac{\epsilon(t)}{2}\|x_{\epsilon(t)}\|^2.
	\end{aligned}
	\]
	Using
	$\mathcal{L}_{\rho}(x^*,\lambda^\infty)
	\leq\mathcal{L}_{\rho}(x_{\epsilon(t)},\lambda^\infty)$
	and taking $x=x(t)+\beta(t)\dot{x}(t)$, we deduce that
	\[
	\begin{aligned}
		\widetilde{\mathcal{E}}(t)
		\geq{}&
		t^{\frac{\alpha+3}{3}}\frac{\epsilon(t)}{2}
		\Big(
		\|x(t)+\beta(t)\dot{x}(t)-x_{\epsilon(t)}\|^2
		+\|x_{\epsilon(t)}\|^2-\|x^*\|^2
		+\|x(t)\|^2-\|x(t)+\beta(t)\dot{x}(t)\|^2
		\Big)\\
		={}&t^{\frac{\alpha+3}{3}}\frac{\epsilon(t)}{2}
		\Big(
		\|x(t)-x_{\epsilon(t)}\|^2
		+\|x_{\epsilon(t)}\|^2-\|x^*\|^2
		\Big)
		-t^{\frac{\alpha+3}{3}}\beta(t)\epsilon(t)
		\langle\dot{x}(t),x_{\epsilon(t)}\rangle.
	\end{aligned}
	\]
	Combining this estimate with~\eqref{DS_thm2_1} yields
	\[
	\begin{aligned}
		\|x(t)-x_{\epsilon(t)}\|^2
		\leq{}&
		\|x^*\|^2-\|x_{\epsilon(t)}\|^2
		+2\beta(t)\langle\dot{x}(t),x_{\epsilon(t)}\rangle\\
		&+\frac{2\widetilde{\mathcal{E}}(t_3)}
		{t^{\frac{\alpha+3}{3}}\epsilon(t)}
		+\frac{2C_6\|x^*\|^2}
		{t^{\frac{\alpha+3}{3}}\epsilon(t)}
		\int_{t_3}^{t}
		s^{\frac{\alpha+3}{3}}\epsilon(s)^2\,ds\\
		&+\frac{2\alpha(\alpha-3)}
		{9t^{\frac{\alpha+3}{3}}\epsilon(t)}
		\int_{t_3}^{t}
		s^{\frac{\alpha-6}{3}}
		\|\lambda(s)-\lambda^\infty\|^2\,ds.
	\end{aligned}
	\]
	By Lemma~\ref{lem:tikhonov_AL},
	$x_{\epsilon(t)}\to x^*$.
	Together with $\dot{x}(t)\to0$ and boundedness of
	$\beta(t)$, this shows that the first two terms
	on the right-hand side tend to zero.
	The next two terms vanish by the assumptions on
	$\epsilon$.
	For the last term, the lower bound on
	$t^2\epsilon(t)$ gives
	\[
	\begin{aligned}
		0\leq{}
		\frac{2\alpha(\alpha-3)}
		{9t^{\frac{\alpha+3}{3}}\epsilon(t)}
		\int_{t_3}^{t}
		s^{\frac{\alpha-6}{3}}
		\|\lambda(s)-\lambda^\infty\|^2\,ds
		\leq{}
		\frac{1}{t^{\frac{\alpha-3}{3}}}
		\int_{t_3}^{t}
		s^{\frac{\alpha-6}{3}}
		\|\lambda(s)-\lambda^\infty\|^2\,ds
		\longrightarrow0.
	\end{aligned}
	\]
	The last limit follows by splitting the integral at a
	fixed sufficiently large time and using
	$\lambda(s)\to\lambda^\infty$ and $\alpha>3$.
	Consequently,
	$\|x(t)-x_{\epsilon(t)}\|\to0$, and the triangle
	inequality yields $x(t)\to x^*$.
	
	\medskip
	\noindent
	\textbf{Case II:}
	There exists $T\geq t_0$ such that
	$\|x(t)\|<\|x^*\|$ for all $t\geq T$.
	
	The trajectory $x(t)$ is bounded.
	Since $\beta(t)\dot{x}(t)\to0$, continuity of
	$\mathcal{L}_{\rho}(\cdot,\lambda^\infty)$ on bounded
	sets and the gap convergence established above imply
	\begin{equation}\label{DS_thm2_3}
		\mathcal{L}_{\rho}(x(t),\lambda^\infty)
		-\mathcal{L}_{\rho}(x^*,\lambda^\infty)\to0.
	\end{equation}
	Let $\tilde{x}$ be any cluster point, and choose
	$t_n\to+\infty$ such that $x(t_n)\to\tilde{x}$.
	By~\eqref{DS_thm2_3} and continuity,
	\[
	\mathcal{L}_{\rho}(\tilde{x},\lambda^\infty)
	=\mathcal{L}_{\rho}(x^*,\lambda^\infty),
	\qquad
	\|\tilde{x}\|\leq\|x^*\|.
	\]
	As shown in the proof of Lemma~\ref{lem:tikhonov_AL},
	$\operatorname*{argmin}\mathcal{L}_{\rho}
	(\cdot,\lambda^\infty)=S$.
	Hence $\tilde{x}\in S$, and the uniqueness of the
	minimum-norm solution gives $\tilde{x}=x^*$.
	Thus the bounded trajectory has the unique cluster
	point $x^*$ and converges to it in $\mathbb R^n$.
	
	\medskip
	\noindent
	\textbf{Case III:}
	Neither Case I nor Case II holds.
	
	For every $T\geq t_0$, there exist $t,s\geq T$
	such that $\|x(t)\|<\|x^*\|$ and
	$\|x(s)\|\geq\|x^*\|$.
	By continuity of the trajectory, there exists a
	sequence $t_n\to+\infty$ such that
	\begin{equation}\label{DS_thm2_4}
		\|x(t_n)\|=\|x^*\|,
		\qquad n\in\mathbb N.
	\end{equation}
	The sequence $\{x(t_n)\}$ is bounded, and
	$\beta(t_n)\dot{x}(t_n)\to0$.
	Applying the argument of Case II along this sequence
	shows that every cluster point belongs to $S$ and
	has norm $\|x^*\|$.
	Therefore, $x(t_n)\to x^*$, and
	\[
	\liminf_{t\to+\infty}\|x(t)-x^*\|=0.
	\]
	This completes the proof.
\end{proof}
\begin{remark}
	Set $t_0=1$ and choose
	\[
	\alpha=4,\qquad M=2,\qquad \rho=1,\qquad
	(\gamma,\beta)=(1,1),\qquad
	\epsilon(t)=\frac{a}{t^2},\quad a\geq1.
	\]
	These choices satisfy all parameter and regularization
	assumptions of Theorem~\ref{DS_thm2}.
	Together with Theorem~\ref{thm:critical_epsilon},
	they yield both $\mathcal{O}(t^{-2})$ convergence rates
	and minimum-norm strong convergence under the trajectory
	condition in Theorem~\ref{DS_thm2}.
\end{remark}

\section{An Implicit Primal--Dual Algorithm and Its
	Convergence Analysis}\label{Algorithm_section}

We derive an implicit primal--dual algorithm from the
dynamical system~\eqref{DS} and analyze its convergence.
Setting $y(t)=x(t)+\frac{3t}{2\alpha}\dot{x}(t)$ and
$z(t)=x(t)+\beta(t)\dot{x}(t)$, we rewrite
the dynamical system~\eqref{DS} as
\[
\begin{cases}
	\displaystyle
	\dot{y}(t)
	=\frac{3-\alpha}{2\alpha}\dot{x}(t)
	-\frac{3t}{2\alpha}
	\big(
	\nabla_x\mathcal{L}_{\rho}(z(t),\lambda(t))
	+\epsilon(t)x(t)
	\big),\\[1mm]
	\displaystyle
	\dot{\lambda}(t)
	=tA\big(y(t)+\beta(t)\dot{y}(t)\big)
	+\frac{3t}{2\alpha}
	\big(t\dot{\beta}(t)-\beta(t)\big)A\dot{x}(t)-tb.
\end{cases}
\]
For the discrete analysis, $f:\mathbb R^n\to\mathbb R$
is convex and possibly nonsmooth, and
$\Omega\neq\emptyset$ is retained with stationarity
$0\in\partial f(x^*)+A^\top\lambda^*$.
On the grid $t_k=kh$, let $s=h^2$ with $h>0$ and denote
the discrete approximations by $(x_k,y_k,z_k,\lambda_k)$.
The backward difference $(x_k-x_{k-1})/h$ gives
$\beta_k=\beta(t_k)/h=\gamma/h+\beta/(ks)$;
set $\epsilon_k=\epsilon(t_k)$.
Introducing an index shift $\theta\geq2$ and replacing
$\nabla f$ by $\partial f$, we obtain the implicit scheme
for integers $k>\theta$:
\begin{equation}\label{discretization_scheme}
	\left\{
	\begin{aligned}
		\omega_k
		&=\frac{(3k+3-3\theta+2\alpha)ks}
		{2\alpha(k-\theta)}
		\big((k-\theta)\beta_{k+1}
		-(k+1-\theta)\beta_k\big),\\
		y_k
		&=x_k+\frac{3(k-\theta)}{2\alpha}(x_k-x_{k-1}),\\
		z_k
		&=x_k+\beta_k(x_k-x_{k-1}),\\
		x_{k+1}
		&\in x_k+\frac{2\alpha}{3-\alpha}(y_{k+1}-y_k)
		+\frac{3ks}{3-\alpha}
		\big(
		\partial f(z_{k+1})+A^\top\lambda_{k+1}
		+\rho A^\top(Az_{k+1}-b)+\epsilon_kx_{k+1}
		\big),\\
		\lambda_{k+1}
		&=\lambda_k
		+\left(1+\frac{k+1-\theta}{k-\theta}\beta_k\right)
		ksAy_{k+1}
		-\frac{k(k+1-\theta)}{k-\theta}\beta_ksAy_k
		-ksb+\omega_kA(x_{k+1}-x_k).
	\end{aligned}
	\right.
\end{equation}
The subsequent analysis allows any regularization sequence
$\{\epsilon_k\}$ satisfying the stated discrete assumptions.
Scheme~\eqref{discretization_scheme} is implemented through
the minimization subproblem in Algorithm~\ref{Algorithm1}.

\begin{algorithm}[H]
	\caption{Implicit primal--dual algorithm}
	\label{Algorithm1}
	\SetAlgoLined
	\DontPrintSemicolon
	
	\textbf{Initialization:}
	Choose $\alpha>3$, $2\leq\theta<3$, $\rho>0$, $s=h^2>0$,
	and either $\gamma>0$, $\beta\in\mathbb R$, or
	$\gamma=0$, $\beta>0$.
	Set $\beta_k=\gamma/h+\beta/(ks)$ and choose an integer
	$K_0>\theta$ such that $1+\beta_k>0$ for all $k\geq K_0$.
	Given $x_0\in\mathbb R^n$ and $\lambda_0\in\mathbb R^m$,
	set $x_{K_0-1}=x_{K_0}=y_{K_0}=z_{K_0}=x_0$
	and $\lambda_{K_0}=\lambda_0$.\;
	
\For{$k=K_0,K_0+1,\ldots$}{
	\textbf{Step 1:}
	Choose $\epsilon_k\geq0$ and compute $\omega_k$
	from~\eqref{discretization_scheme}.
	Set $\theta_k=\rho+\frac{2\alpha+3(k+1-\theta)}{2\alpha}ks$ and
\[
\begin{aligned}
	&\delta_k=\frac{k+\alpha-\theta+ks\epsilon_k}{(1+\beta_{k+1})ks},
	\quad \sigma_k=ks+\rho-\theta_k,
	\quad \zeta_k=\epsilon_k-\delta_k,\\
	&e_k=\frac{\sigma_k}{\theta_k}Ax_k
	-\frac{3k(k+1-\theta)}{2\alpha\theta_k}
	\beta_ksA(x_k-x_{k-1}).
\end{aligned}
\]
	Compute
	\[
	z_{k+1}=\operatorname*{argmin}_{z\in\mathbb R^n}
	\left\{
	f(z)
	+\frac{\theta_k}{2}
	\left\|Az+e_k-\frac{ks+\rho}{\theta_k}b
	+\frac{\lambda_k}{\theta_k}\right\|^2
	+\frac{\delta_k}{2}
	\left\|z+\frac{\zeta_k}{\delta_k}x_k
	-\frac{k-\theta}{ks\delta_k}(x_k-x_{k-1})\right\|^2
	\right\}.
	\]
	
	\textbf{Step 2:}
	Set $x_{k+1}=(z_{k+1}+\beta_{k+1}x_k)/(1+\beta_{k+1})$.\;
	
	\textbf{Step 3:}
	Set $y_{k+1}=x_{k+1}
	+\frac{3(k+1-\theta)}{2\alpha}(x_{k+1}-x_k)$ and update
	\[
	\lambda_{k+1}=\lambda_k
	+\left(1+\frac{k+1-\theta}{k-\theta}\beta_k\right)ksAy_{k+1}
	-\frac{k(k+1-\theta)}{k-\theta}\beta_ksAy_k
	-ksb+\omega_kA(x_{k+1}-x_k).
	\]
	\If{a stopping criterion is satisfied}{
		\Return $(x_{k+1},\lambda_{k+1})$.\;
	}
}
\end{algorithm}

The initialization ensures $\theta_k>0$ and $\delta_k>0$,
so the subproblem in Step~1 is strongly convex and has
a unique minimizer.

\begin{proposition}
	Algorithm~\ref{Algorithm1} is equivalent to
	scheme~\eqref{discretization_scheme} with the same
	parameters and initial values.
\end{proposition}

\begin{proof}
	The optimality condition for Step~1 of
	Algorithm~\ref{Algorithm1} is
	\[
	\begin{aligned}
		0\in{}\partial f(z_{k+1})
		+\theta_kA^\top
		\left(
		Az_{k+1}+e_k-\frac{ks+\rho}{\theta_k}b
		+\frac{\lambda_k}{\theta_k}
		\right)
		+\delta_k
		\left(
		z_{k+1}+\frac{\zeta_k}{\delta_k}x_k
		-\frac{k-\theta}{ks\delta_k}(x_k-x_{k-1})
		\right).
	\end{aligned}
	\]
	Substituting the definitions of $\theta_k$ and $e_k$
	and using Steps~2--3, we obtain
	\[
	\begin{aligned}
		0\in{}\partial f(z_{k+1})
		+A^\top\lambda_{k+1}
		+\rho A^\top(Az_{k+1}-b)
		+\delta_kz_{k+1}+\zeta_kx_k
		-\frac{k-\theta}{ks}(x_k-x_{k-1}).
	\end{aligned}
	\]
	Using $\zeta_k=\epsilon_k-\delta_k$,
	$z_{k+1}=(1+\beta_{k+1})x_{k+1}-\beta_{k+1}x_k$,
	and the definitions of $\delta_k$, $y_k$, and $y_{k+1}$
	gives
	\[
	\begin{aligned}
		0\in{}
		\frac{3-\alpha}{3ks}(x_k-x_{k+1})
		+\frac{2\alpha}{3ks}(y_{k+1}-y_k)
		+\partial f(z_{k+1})
		+A^\top\lambda_{k+1}
		+\rho A^\top(Az_{k+1}-b)
		+\epsilon_kx_{k+1}.
	\end{aligned}
	\]
	Rearranging yields the primal inclusion
	in~\eqref{discretization_scheme}; the remaining relations
	follow directly from the algorithm.
	Conversely, these substitutions are reversible, and
	the resulting optimality condition characterizes the
	unique minimizer in Step~1.
	Hence the two formulations are equivalent.
\end{proof}

\begin{lemma}\label{AL_lem1}
	Let $\{(x_k,\lambda_k)\}_{k\geq K_0}$ be generated by
	Algorithm~\ref{Algorithm1}, and let
	$(x^*,\kappa)\in\mathbb F\times\mathbb R^m$.
	Define
	\begin{equation}\label{AL_lem1_1}
		\begin{aligned}
			\mathcal E_\kappa^k
			:={}&k(k+1-\theta)s
			\left(
			\mathcal L_\rho(z_k,\kappa)
			-\mathcal L_\rho(x^*,\kappa)
			+\frac{\epsilon_k}{2}\|x_k\|^2
			\right)\\
			&+\frac12\left\|\frac{2\alpha}{3}(y_k-x^*)\right\|^2
			+\frac{\alpha(\alpha-3)}{9}\|x_k-x^*\|^2
			+\frac{\alpha}{3}\|\lambda_k-\kappa\|^2.
		\end{aligned}
	\end{equation}
	Suppose that $\alpha>3$, $\theta\geq2$, and
	$\{\epsilon_k\}_{k\geq K_0}$ is nonnegative and
	nonincreasing. Assume that there exist $M>1$ and
	an integer $K_1\geq K_0$ such that
	\begin{equation}\label{AL_parameter_condition_1}
		2(k+1)(k+2-\theta)(\epsilon_{k+1}-\epsilon_k)
		+\frac{Mk^2(k+1-\theta)\epsilon_k^2|\beta_k|s}
		{k-\theta}\leq0,
		\qquad k\geq K_1.
	\end{equation}
	Then there exist $K_2\geq K_1$ and $r_1>0$ such that,
	for all $k\geq K_2$,
	\begin{equation}\label{AL_lem1_conclusion}
		\begin{aligned}
			\mathcal E_\kappa^{k+1}-\mathcal E_\kappa^k
			\leq{}&
			\left(\frac{2(3-\alpha)}{3}k+2-\theta\right)s
			\big(
			\mathcal L_\rho(z_{k+1},\kappa)
			-\mathcal L_\rho(x^*,\kappa)
			\big)\\
			&-r_1(k+1-\theta)\|x_{k+1}-x_k\|^2
			+\frac{\alpha}{3}ks\epsilon_k
			\big(\|x^*\|^2-\|x_{k+1}-x^*\|^2\big)\\
			&+\left(\frac{3-\alpha}{3}k+1-\frac{\theta}{2}\right)
			s\epsilon_k\|x_{k+1}\|^2.
		\end{aligned}
	\end{equation}
\end{lemma}

\begin{proof}
	By~\eqref{discretization_scheme}, there exists
	$\xi_{k+1}\in\partial f(z_{k+1})$ such that
	\begin{equation}\label{mid_1}
		\begin{aligned}
			&-\frac{2\alpha}{3ks}(y_{k+1}-y_k)
			+\frac{3-\alpha}{3ks}(x_{k+1}-x_k)
			-\epsilon_kx_{k+1}
			-A^\top(\lambda_{k+1}-\kappa)\\
			&\qquad=\xi_{k+1}+A^\top\kappa
			+\rho A^\top(Az_{k+1}-b).
		\end{aligned}
	\end{equation}
	Using the definition of $y_k$, we obtain
	\begin{equation}\label{mid_2}
		\begin{aligned}
			x_k-x_{k-1}
			={}\left(1+\frac{\alpha}{k-\theta}\right)
			(x_{k+1}-x_k)
			+\frac{ks}{k-\theta}
			\big(
			\xi_{k+1}+A^\top\lambda_{k+1}
			+\rho A^\top(Az_{k+1}-b)+\epsilon_kx_{k+1}
			\big).
		\end{aligned}
	\end{equation}
	
	The subgradient inequality and the quadratic penalty give
	\begin{equation}\label{mid_3}
		\begin{aligned}
			\mathcal L_\rho(z_{k+1},\kappa)
			-\mathcal L_\rho(z_k,\kappa)
			\leq{}
			\left\langle
			\xi_{k+1}+A^\top\kappa+\rho A^\top(Az_{k+1}-b),
			z_{k+1}-z_k
			\right\rangle
			-\frac{\rho}{2}\|A(z_{k+1}-z_k)\|^2.
		\end{aligned}
	\end{equation}
	Similarly, using $Ax^*=b$ and
	$z_{k+1}=x_{k+1}+\beta_{k+1}(x_{k+1}-x_k)$ yields
	\begin{equation}\label{mid_4}
		\begin{aligned}
			&\left\langle
			\xi_{k+1}+A^\top\kappa+\rho A^\top(Az_{k+1}-b)
			+\epsilon_kx_{k+1},x^*-x_{k+1}
			\right\rangle\\
			&\quad\leq
			\mathcal L_\rho(x^*,\kappa)
			-\mathcal L_\rho(z_{k+1},\kappa)
			-\frac{\rho}{2}\|Az_{k+1}-b\|^2\\
			&\qquad
			+\frac{\epsilon_k}{2}
			\big(
			\|x^*\|^2-\|x_{k+1}\|^2-\|x_{k+1}-x^*\|^2
			\big)\\
			&\qquad
			+\beta_{k+1}
			\left\langle
			\xi_{k+1}+A^\top\lambda_{k+1}
			+\rho A^\top(Az_{k+1}-b),x_{k+1}-x_k
			\right\rangle\\
			&\qquad
			+\beta_{k+1}
			\langle\kappa-\lambda_{k+1},A(x_{k+1}-x_k)\rangle.
		\end{aligned}
	\end{equation}
	
	We now expand the squared-norm differences in
	$\mathcal E_\kappa^{k+1}-\mathcal E_\kappa^k$ using
	$\frac12\|u\|^2-\frac12\|v\|^2
	=\langle u,u-v\rangle-\frac12\|u-v\|^2$.
	Substitute~\eqref{mid_1} and the dual update, and use
	\eqref{mid_2} in~\eqref{mid_3} together with
	$z_{k+1}-z_k=(1+\beta_{k+1})(x_{k+1}-x_k)
	-\beta_k(x_k-x_{k-1})$.
	Combining the resulting estimates with~\eqref{mid_4},
	the mixed terms involving $\lambda_{k+1}-\kappa$ cancel.
	After discarding nonpositive squared-norm terms,
	we obtain
	\[
	\begin{aligned}
		\mathcal E_\kappa^{k+1}-\mathcal E_\kappa^k
		\leq{}&
		\left(\frac{2(3-\alpha)}{3}k+2-\theta\right)s
		\big(
		\mathcal L_\rho(z_{k+1},\kappa)
		-\mathcal L_\rho(x^*,\kappa)
		\big)\\
		&+\mathcal C_k^1\|x_{k+1}-x_k\|^2
		-\frac{\alpha\rho}{3}ks\|Az_{k+1}-b\|^2
		+\mathcal C_k^2\|x_{k+1}\|^2\\
		&+\frac{\alpha}{3}ks\epsilon_k
		\big(\|x^*\|^2-\|x_{k+1}-x^*\|^2\big)\\
		&-\frac{k^2(k+1-\theta)}{k-\theta}\beta_ks^2
		\|\xi_{k+1}+A^\top\lambda_{k+1}
		+\rho A^\top(Az_{k+1}-b)\|^2\\
		&-\frac{k^2(k+1-\theta)}{k-\theta}\beta_ks^2\epsilon_k
		\left\langle
		\xi_{k+1}+A^\top\lambda_{k+1}
		+\rho A^\top(Az_{k+1}-b),x_{k+1}
		\right\rangle\\
		&+\mathcal C_k^3
		\left\langle
		\xi_{k+1}+A^\top\lambda_{k+1}
		+\rho A^\top(Az_{k+1}-b),x_{k+1}-x_k
		\right\rangle,
	\end{aligned}
	\]
	where
	\[
	\begin{aligned}
		\mathcal C_k^1
		&:=\frac{3-\alpha}{3}(k+1-\theta)
		-\frac{\alpha(\alpha-3)}{9}
		-\frac{k(k+1-\theta)s\epsilon_k}{2},\\
		\mathcal C_k^2
		&:=\frac{(k+1)(k+2-\theta)s}{2}
		(\epsilon_{k+1}-\epsilon_k)
		+\left(\frac{3-\alpha}{3}k+1-\frac{\theta}{2}\right)
		s\epsilon_k,\\
		\mathcal C_k^3
		&:=\left(
		\frac{2\alpha}{3}
		\left(1-\frac{3(k+1-\theta)}{2(k-\theta)}\right)\beta_{k+1}
		+\frac{(k+1-\theta)(k+\alpha-\theta)}{k-\theta}
		(\beta_{k+1}-\beta_k)
		\right)ks.
	\end{aligned}
	\]
	
	Apply Young's inequality to the last two inner products,
	using $M$ from~\eqref{AL_parameter_condition_1} for the
	regularization term and the weights $\sqrt{k}s$
	and $1/(\sqrt{k}s)$ for the term with $\mathcal C_k^3$.
	Since
	$\mathcal C_k^1\leq(3-\alpha)(k+1-\theta)/3$,
	this gives
	\begin{equation}\label{AL_lem1_7}
		\begin{aligned}
			\mathcal E_\kappa^{k+1}-\mathcal E_\kappa^k
			\leq{}&
			\left(\frac{2(3-\alpha)}{3}k+2-\theta\right)s
			\big(
			\mathcal L_\rho(z_{k+1},\kappa)
			-\mathcal L_\rho(x^*,\kappa)
			\big)\\
			&+\left(
			\frac{3-\alpha}{3}(k+1-\theta)
			+\frac{|\mathcal C_k^3|}{2\sqrt{k}s}
			\right)\|x_{k+1}-x_k\|^2\\
			&-\frac{\alpha\rho}{3}ks\|Az_{k+1}-b\|^2
			+\frac{\alpha}{3}ks\epsilon_k
			\big(\|x^*\|^2-\|x_{k+1}-x^*\|^2\big)\\
			&+\left(
			\mathcal C_k^2
			+\frac{Mk^2(k+1-\theta)\epsilon_k^2|\beta_k|s^2}
			{4(k-\theta)}
			\right)\|x_{k+1}\|^2\\
			&-\left(
			\left(\beta_k-\frac{|\beta_k|}{M}\right)
			\frac{k^2(k+1-\theta)s^2}{k-\theta}
			-\frac{\sqrt{k}s}{2}|\mathcal C_k^3|
			\right)\\
			&\qquad\cdot
			\|\xi_{k+1}+A^\top\lambda_{k+1}
			+\rho A^\top(Az_{k+1}-b)\|^2.
		\end{aligned}
	\end{equation}
	
	For $\beta_k=\gamma/h+\beta/(ks)$, we have
	$\mathcal C_k^3=\mathcal O(k)$ if $\gamma>0$ and
	$\mathcal C_k^3=\mathcal O(1)$ if $\gamma=0$, $\beta>0$.
	Thus, since $\alpha>3$ and $M>1$, there exist
	$K_1'\geq K_0$ and $r_1>0$ such that, for $k\geq K_1'$,
	\[
	\begin{aligned}
		\frac{3-\alpha}{3}(k+1-\theta)
		+\frac{|\mathcal C_k^3|}{2\sqrt{k}s}
		&\leq-r_1(k+1-\theta),\\
		\left(\beta_k-\frac{|\beta_k|}{M}\right)
		\frac{k^2(k+1-\theta)s^2}{k-\theta}
		-\frac{\sqrt{k}s}{2}|\mathcal C_k^3|
		&\geq0.
	\end{aligned}
	\]
	Moreover, condition~\eqref{AL_parameter_condition_1}
	implies
	\[
	\mathcal C_k^2
	+\frac{Mk^2(k+1-\theta)\epsilon_k^2|\beta_k|s^2}
	{4(k-\theta)}
	\leq
	\left(\frac{3-\alpha}{3}k+1-\frac{\theta}{2}\right)
	s\epsilon_k,
	\qquad k\geq K_1.
	\]
	Substituting these estimates into~\eqref{AL_lem1_7}
	and discarding the nonpositive penalty term proves
	\eqref{AL_lem1_conclusion} for
	$K_2=\max\{K_1,K_1'\}$.
\end{proof}

The following lemma can be viewed as a discrete counterpart of Lemma \ref{lemB.2}, and it will be instrumental in establishing the convergence rate and the strong convergence of the sequence generated by our proposed discrete scheme.
\begin{lemma}(Kronecker's Lemma){\rm{\cite[Lemma IV.3.2]{Shiryaev Springer 1996}}}\label{lemB.4}
	Suppose that $(\phi_k)_{k=K_0}^{\infty}$ is an infinite sequence of real numbers such that $\sum_{k=K_0}^{\infty}\phi_k=c$ exists and is finite, where $K_0\ge 1$. And let $(\psi_k)_{k\ge K_0}$ be a nondecreasing function in $\left(0,+\infty\right)$ such that $\lim_{k \to +\infty}\psi_k=+\infty$. Then,
	\begin{equation*}
		\lim_{k \to +\infty}\frac{1}{\psi_k}\sum_{j=K_0}^{k}\psi_j\phi_j=0.
	\end{equation*}
\end{lemma}

\begin{theorem}\label{AL_thm1}
	Suppose that the assumptions of Lemma~\ref{AL_lem1} hold,
	$2\leq\theta<3$, and
	$\sum_{k=K_0}^{+\infty}\epsilon_k/k<+\infty$.
	Let $\{(x_k,\lambda_k)\}_{k\geq K_0}$ be generated by
	Algorithm~\ref{Algorithm1}, and let
	$(x^*,\lambda^*)\in\Omega$.
	Then the following statements hold.
	
	\noindent
	(i) As $k\to+\infty$,
	\[
	\mathcal L_\rho(z_k,\lambda^*)
	-\mathcal L_\rho(x^*,\lambda^*)\to0,
	\qquad
	\|x_{k+1}-x_k\|\to0,
	\qquad
	\|Az_k-b\|\to0.
	\]
	
	\noindent
	(ii) If, in addition,
	$\sum_{k=K_0}^{+\infty}k\epsilon_k<+\infty$,
	then $\{x_k\}$ is bounded and
	\[
	\begin{aligned}
		\mathcal L_\rho(z_{k+1},\lambda^*)
		-\mathcal L_\rho(x^*,\lambda^*)
		&=\mathcal O(k^{-2}),
		&\qquad \|x_{k+1}-x_k\|&=\mathcal O(k^{-1}),\\
		\|Az_{k+1}-b\|
		&=\mathcal O(k^{-2}),
		&|f(z_{k+1})-f(x^*)|&=\mathcal O(k^{-2}).
	\end{aligned}
	\]
\end{theorem}

\begin{proof}
	The saddle-point property gives
	$\mathcal L_\rho(z_k,\lambda^*)
	-\mathcal L_\rho(x^*,\lambda^*)
	\geq\frac{\rho}{2}\|Az_k-b\|^2\geq0$.
	Summing~\eqref{AL_lem1_conclusion} with
	$\kappa=\lambda^*$ from $K_2$ to $k-1$ therefore yields
	\begin{equation}\label{AL_thm3.1_ex1}
		0\leq\mathcal E_{\lambda^*}^k
		\leq\mathcal E_{\lambda^*}^{K_2}
		+\frac{\alpha}{3}\|x^*\|^2
		\sum_{i=K_2}^{k-1}is\epsilon_i.
	\end{equation}
	By Lemma~\ref{lemB.4}, with
	$\phi_i=\epsilon_i/i$ and $\psi_i=i^2s$,
	\[
	\frac{1}{k^2s}\sum_{i=K_2}^{k-1}is\epsilon_i\to0,
	\qquad
	\frac{\mathcal E_{\lambda^*}^k}{k(k+1-\theta)s}\to0.
	\]
	The definition of $\mathcal E_{\lambda^*}^k$ implies
	that the augmented-Lagrangian gap tends to zero,
	$(x_k-x^*)/k\to0$, and $(y_k-x^*)/k\to0$.
	Using $x_k-x_{k-1}
	=\frac{2\alpha}{3(k-\theta)}(y_k-x_k)$ and the
	gap lower bound proves (i).
	
	Now assume $\sum_{k=K_0}^{+\infty}k\epsilon_k<+\infty$.
	Estimate~\eqref{AL_thm3.1_ex1} shows that
	$\mathcal E_{\lambda^*}^k$ is bounded.
	Hence $\{x_k\}$, $\{y_k\}$, and $\{\lambda_k\}$
	are bounded, the augmented-Lagrangian gap is
	$\mathcal O(k^{-2})$, and the same identity for
	$x_k-x_{k-1}$ gives
	$\|x_{k+1}-x_k\|=\mathcal O(k^{-1})$.
	
	To obtain the remaining rates, Lemma~\ref{AL_lem1}
	and the saddle-point property give, for any
	$\kappa\in\mathbb R^m$,
	\begin{equation}\label{AL_thm1_1}
		\mathcal E_\kappa^{k+1}-\mathcal E_\kappa^k
		\leq
		\left(\frac{2(3-\alpha)}{3}k+2-\theta\right)s
		\langle\kappa-\lambda^*,Az_{k+1}-b\rangle
		+\frac{\alpha}{3}ks\epsilon_k\|x^*\|^2.
	\end{equation}
	Fix an integer $\tau\geq K_2$ and set
	\[
	\kappa_\tau:=
	\begin{cases}
		\displaystyle
		\lambda^*+\frac{Az_{\tau+1}-b}{\|Az_{\tau+1}-b\|},
		&Az_{\tau+1}-b\neq0,\\
		\lambda^*,&Az_{\tau+1}-b=0.
	\end{cases}
	\]
	Then $\|\kappa_\tau-\lambda^*\|\leq1$.
	Let
	$a_k=a_1(2+\frac{2\alpha}{3}-\theta)_{k-1}/
	[k(3-\theta)_{k-1}]$, where $a_1>0$ and
	$(u)_j$ denotes the rising factorial, with $(u)_0=1$.
	Since $\alpha>3$ and $2\leq\theta<3$, this sequence
	is positive and increasing, and
	\[
	(a_{k+1}-a_k)(k+1)(k+2-\theta)
	=a_k\left(\frac{2(\alpha-3)}{3}k-2+\theta\right).
	\]
	Multiplying~\eqref{AL_thm1_1} by $a_k$ and summing
	by parts therefore yields
	\begin{equation}\label{AL_thm1_2}
		\begin{aligned}
			a_{k+1}\mathcal E_{\kappa_\tau}^{k+1}
			\leq{}&
			a_{K_2}\mathcal E_{\kappa_\tau}^{K_2}
			+\sum_{i=K_2}^{k}(a_{i+1}-a_i)
			\Big(
			\mathcal E_{\kappa_\tau}^{i+1}
			-(i+1)(i+2-\theta)s
			\langle\kappa_\tau-\lambda^*,Az_{i+1}-b\rangle
			\Big)\\
			&+\frac{\alpha}{3}\|x^*\|^2
			\sum_{i=K_2}^{k}a_iis\epsilon_i.
		\end{aligned}
	\end{equation}
	
	By boundedness of $\mathcal E_{\lambda^*}^k$ and
	$\|\lambda_k-\kappa_\tau\|^2
	\leq2\|\lambda_k-\lambda^*\|^2+2$, there exists
	$\widetilde C>0$, independent of $k$ and $\tau$, such that
	\[
	\begin{aligned}
		\mathcal E_{\kappa_\tau}^k
		&\leq2\mathcal E_{\lambda^*}^k
		+k(k+1-\theta)s
		\langle\kappa_\tau-\lambda^*,Az_k-b\rangle
		+\frac{2\alpha}{3}\\
		&\leq\widetilde C+k(k+1-\theta)s
		\langle\kappa_\tau-\lambda^*,Az_k-b\rangle.
	\end{aligned}
	\]
	In particular,
	$\mathcal E_{\kappa_\tau}^{K_2}
	\leq\widetilde C+K_2(K_2+1-\theta)s\|Az_{K_2}-b\|
	=:\widetilde C_1$ uniformly in $\tau$.
	Substituting the preceding bound into~\eqref{AL_thm1_2}
	gives
	\[
	\begin{aligned}
		\mathcal E_{\kappa_\tau}^{k+1}
		\leq{}&
		\frac{a_{K_2}}{a_{k+1}}
		\mathcal E_{\kappa_\tau}^{K_2}
		+\left(1-\frac{a_{K_2}}{a_{k+1}}\right)\widetilde C
		+\frac{\alpha\|x^*\|^2}{3a_{k+1}}
		\sum_{i=K_2}^{k}a_iis\epsilon_i\\
		\leq{}&
		\widetilde C_1+\widetilde C
		+\frac{\alpha}{3}\|x^*\|^2
		\sum_{i=K_2}^{+\infty}is\epsilon_i
		=:\widetilde C_4<+\infty,
	\end{aligned}
	\]
	where we used $0<a_i\leq a_{k+1}$ for $i\leq k$.
	The constant $\widetilde C_4$ is independent of both
	$k$ and $\tau$.
	
	Taking $\tau=k$ and dropping the nonnegative terms
	in $\mathcal E_{\kappa_k}^{k+1}$, we obtain
	\[
	(k+1)(k+2-\theta)s
	\big(
	\mathcal L(z_{k+1},\lambda^*)
	-\mathcal L(x^*,\lambda^*)
	+\|Az_{k+1}-b\|
	\big)
	\leq\widetilde C_4,
	\]
	where $\mathcal L(x,\lambda)=f(x)+\langle\lambda,Ax-b\rangle$.
	Since
	$\mathcal L(z_{k+1},\lambda^*)-\mathcal L(x^*,\lambda^*)
	\geq0$, it follows that
	\[
	\|Az_{k+1}-b\|
	\leq\frac{\widetilde C_4}{(k+1)(k+2-\theta)s}.
	\]
	The same estimate and the saddle-point property give
	\[
	-\|\lambda^*\|\|Az_{k+1}-b\|
	\leq f(z_{k+1})-f(x^*)
	\leq
	\frac{\widetilde C_4}{(k+1)(k+2-\theta)s}
	+\|\lambda^*\|\|Az_{k+1}-b\|.
	\]
	Thus both the feasibility violation and the absolute
	objective residual are $\mathcal O(k^{-2})$,
	completing the proof.
\end{proof}

\begin{lemma}\label{AL_conv_dual}
	Suppose that the assumptions of Lemma~\ref{AL_lem1} hold,
	$2\leq\theta<3$, and
	$\sum_{k=K_0}^{+\infty}\epsilon_k/k<+\infty$.
	Assume additionally that $f$ is continuously differentiable
	with an $L$-Lipschitz continuous gradient.
	Let $\{(x_k,\lambda_k)\}_{k\geq K_0}$ be generated by
	Algorithm~\ref{Algorithm1}.
	Then there exists $(x^*,\lambda^\infty)\in\Omega$
	such that $\lambda_k\to\lambda^\infty$.
\end{lemma}

\begin{proof}
	If $\operatorname{Ran}A=\{0\}$, feasibility implies $b=0$,
	and the dual update gives $\lambda_{k+1}=\lambda_k$.
	The conclusion then follows immediately.
	Assume henceforth that $\operatorname{Ran}A\neq\{0\}$.
	
	Take $(x^*,\bar\lambda)\in\Omega$ and set
	$\lambda^\infty
	=P_{\operatorname{Ran}A}\bar\lambda
	+P_{\ker A^\top}\lambda_{K_0}$.
	Then $A^\top\lambda^\infty=A^\top\bar\lambda$, so
	$(x^*,\lambda^\infty)\in\Omega$.
	Since $\lambda_{k+1}-\lambda_k\in\operatorname{Ran}A$,
	the component in $\ker A^\top$ remains constant, and
	$\lambda_k-\lambda^\infty\in\operatorname{Ran}A$.
	
	By Theorem~\ref{AL_thm1}(i),
	\[
	x_{k+1}-x_k\to0,\qquad
	Az_k-b\to0,\qquad
	\mathcal L_\rho(z_k,\bar\lambda)
	-\mathcal L_\rho(x^*,\bar\lambda)\to0.
	\]
	Boundedness of $\{\beta_k\}$ gives
	$z_k-x_k=\beta_k(x_k-x_{k-1})\to0$, and hence
	$Ax_k-b\to0$.
	Moreover, convexity and Lipschitz continuity of $\nabla f$,
	together with $\nabla f(x^*)=-A^\top\bar\lambda$, imply
	\[
	\begin{aligned}
		\|\nabla f(z_k)-\nabla f(x^*)\|^2
		&\leq2L\big(
		f(z_k)-f(x^*)
		-\langle\nabla f(x^*),z_k-x^*\rangle
		\big)\\
		&\leq2L\big(
		\mathcal L_\rho(z_k,\bar\lambda)
		-\mathcal L_\rho(x^*,\bar\lambda)
		\big)\to0.
	\end{aligned}
	\]
	
	Applying $A$ to~\eqref{mid_2} and using
	$\nabla f(x^*)+A^\top\lambda^\infty=0$, we obtain
	\[
	\begin{aligned}
		AA^\top(\lambda_{k+1}-\lambda^\infty)
		={}&\frac{k-\theta}{ks}A(x_k-x_{k-1})
		-\frac{k+\alpha-\theta}{ks}A(x_{k+1}-x_k)\\
		&-A\big(\nabla f(z_{k+1})-\nabla f(x^*)\big)
		-\rho AA^\top(Az_{k+1}-b)
		-\epsilon_kAx_{k+1}.
	\end{aligned}
	\]
	The first four terms on the right-hand side tend to zero
	by the preceding limits.
	Since $\{\epsilon_k\}$ is nonnegative and nonincreasing,
	$\sum_{k=K_0}^{+\infty}\epsilon_k/k<+\infty$
	implies $\epsilon_k\to0$.
	Together with $Ax_{k+1}\to b$, this gives
	$\epsilon_kAx_{k+1}\to0$.
	Thus $AA^\top(\lambda_{k+1}-\lambda^\infty)\to0$.
	
	Finally, $AA^\top$ is positive definite on
	$\operatorname{Ran}A$, so there exists $\sigma_A>0$
	such that
	$\|AA^\top u\|\geq\sigma_A\|u\|$
	for all $u\in\operatorname{Ran}A$.
	As $\lambda_{k+1}-\lambda^\infty$ belongs to this
	subspace, we conclude that
	\[
	\|\lambda_{k+1}-\lambda^\infty\|
	\leq\frac{1}{\sigma_A}
	\|AA^\top(\lambda_{k+1}-\lambda^\infty)\|
	\to0.
	\]
\end{proof}

\begin{theorem}
	\label{thm:discrete_critical_epsilon}
	Suppose that the assumptions of Lemma~\ref{AL_lem1} hold,
	$2\leq\theta<3$, and $\epsilon_k=a/k^2$ with $a>0$.
	Assume additionally that $f$ is continuously differentiable
	with an $L$-Lipschitz continuous gradient.
	Let $\{(x_k,\lambda_k)\}_{k\geq K_0}$ be generated by
	Algorithm~\ref{Algorithm1}.
	Then $\{\lambda_k\}$ is bounded, and for any
	$(x^*,\lambda^*)\in\Omega$,
	\[
	\begin{aligned}
		\mathcal L_\rho(z_{k+1},\lambda^*)
		-\mathcal L_\rho(x^*,\lambda^*)
		&=\mathcal O(k^{-2}),
		&\qquad \|x_{k+1}-x_k\|&=\mathcal O(k^{-1}),\\
		\|Az_{k+1}-b\|
		&=\mathcal O(k^{-2}),
		&|f(z_{k+1})-f(x^*)|&=\mathcal O(k^{-2}).
	\end{aligned}
	\]
\end{theorem}

\begin{proof}
	For $\epsilon_k=a/k^2$, we have
	$\sum_{k=K_0}^{+\infty}\epsilon_k/k<+\infty$.
	Moreover, boundedness of $\{\beta_k\}$ shows that the
	left-hand side of~\eqref{AL_parameter_condition_1}
	is $-4a/k+\mathcal O(k^{-2})$, so this condition holds
	for all sufficiently large $k$.
	Lemma~\ref{AL_conv_dual} therefore gives convergence,
	and hence boundedness, of $\{\lambda_k\}$.
	
	Fix $(x^*,\lambda^*)\in\Omega$.
	Using $y_{k+1}=x_{k+1}
	+\frac{3(k+1-\theta)}{2\alpha}(x_{k+1}-x_k)$
	and $\frac12\|u+v\|^2\leq\|u\|^2+\|v\|^2$,
	we obtain
	\[
	\begin{aligned}
		\frac{c}{k+1}\mathcal E_{\lambda^*}^{k+1}
		\leq{}&
		c(k+2-\theta)s
		\big(
		\mathcal L_\rho(z_{k+1},\lambda^*)
		-\mathcal L_\rho(x^*,\lambda^*)
		\big)\\
		&+\frac{ca(k+2-\theta)s}{2(k+1)^2}\|x_{k+1}\|^2
		+\frac{c\alpha(5\alpha-3)}{9(k+1)}
		\|x_{k+1}-x^*\|^2\\
		&+\frac{c(k+1-\theta)^2}{k+1}
		\|x_{k+1}-x_k\|^2
		+\frac{c\alpha}{3(k+1)}
		\|\lambda_{k+1}-\lambda^*\|^2.
	\end{aligned}
	\]
	Choose
	$0<c\leq\min\{2(\alpha-3)/3,r_1,3as/(5\alpha-3)\}$,
	where $r_1$ is given by Lemma~\ref{AL_lem1}.
	Combining the preceding estimate with
	\eqref{AL_lem1_conclusion}, and using $\theta\geq2$
	and nonnegativity of the augmented-Lagrangian gap,
	gives, for all sufficiently large $k$,
	\[
	\begin{aligned}
		\mathcal E_{\lambda^*}^{k+1}
		-\mathcal E_{\lambda^*}^{k}
		+\frac{c}{k+1}\mathcal E_{\lambda^*}^{k+1}
		\leq{}&
		\left(\frac{2(3-\alpha)}{3}+c\right)(k+2-\theta)s
		\big(
		\mathcal L_\rho(z_{k+1},\lambda^*)
		-\mathcal L_\rho(x^*,\lambda^*)
		\big)\\
		&+\left[
		\left(\frac{3-\alpha}{3}k+1-\frac{\theta}{2}\right)
		\frac{as}{k^2}
		+\frac{ca(k+2-\theta)s}{2(k+1)^2}
		\right]\|x_{k+1}\|^2\\
		&+\left[
		-\frac{\alpha as}{3k}
		+\frac{c\alpha(5\alpha-3)}{9(k+1)}
		\right]\|x_{k+1}-x^*\|^2\\
		&+\left[
		\frac{c(k+1-\theta)^2}{k+1}
		-r_1(k+1-\theta)
		\right]\|x_{k+1}-x_k\|^2\\
		&+\frac{c\alpha}{3(k+1)}
		\|\lambda_{k+1}-\lambda^*\|^2
		+\frac{\alpha as}{3k}\|x^*\|^2.
	\end{aligned}
	\]
	By the choice of $c$, the first four terms on the
	right-hand side are nonpositive.
	Boundedness of $\{\lambda_k\}$ therefore yields
	constants $C>0$ and $K\geq K_2$ such that
	\[
	\mathcal E_{\lambda^*}^{k+1}
	-\mathcal E_{\lambda^*}^{k}
	+\frac{c}{k+1}\mathcal E_{\lambda^*}^{k+1}
	\leq\frac{C}{k+1},
	\qquad k\geq K.
	\]
	Let $P_k=\prod_{j=K}^{k}(1+c/j)$.
	Multiplying this inequality by $P_k$, summing,
	and using $P_{k+1}-P_k=cP_k/(k+1)$, we obtain
	\[
	\mathcal E_{\lambda^*}^{k+1}
	\leq
	\frac{P_K}{P_{k+1}}\mathcal E_{\lambda^*}^{K}
	+\frac{C}{c}
	\left(1-\frac{P_K}{P_{k+1}}\right).
	\]
	Thus $\{\mathcal E_{\lambda^*}^k\}$ is bounded.
	Its definition implies boundedness of $\{x_k\}$
	and $\{y_k\}$ and the estimate
	$\mathcal L_\rho(z_{k+1},\lambda^*)
	-\mathcal L_\rho(x^*,\lambda^*)=\mathcal O(k^{-2})$.
	The identity
	$x_{k+1}-x_k
	=\frac{2\alpha}{3(k+1-\theta)}(y_{k+1}-x_{k+1})$
	also gives $\|x_{k+1}-x_k\|=\mathcal O(k^{-1})$.
	
	For the feasibility and objective residuals, use
	the weights $a_k$ and the vectors $\kappa_\tau$
	introduced in the proof of Theorem~\ref{AL_thm1}.
	The weight identity there implies, for $i\geq K$,
	\[
	\frac{a_i}{i}
	\leq\frac{3}{\alpha-3}(a_{i+1}-a_i).
	\]
	Consequently, the critical regularization satisfies
	\[
	\frac{1}{a_{k+1}}\sum_{i=K}^{k}a_iis\epsilon_i
	=\frac{as}{a_{k+1}}\sum_{i=K}^{k}\frac{a_i}{i}
	\leq\frac{3as}{\alpha-3}
	\left(1-\frac{a_K}{a_{k+1}}\right)
	\leq\frac{3as}{\alpha-3}.
	\]
	Boundedness of $\mathcal E_{\lambda^*}^k$ and
	$\|\kappa_\tau-\lambda^*\|\leq1$ give a constant
	$C_1>0$, independent of $k$ and $\tau$, such that
	\[
	\mathcal E_{\kappa_\tau}^k
	\leq C_1+k(k+1-\theta)s
	\langle\kappa_\tau-\lambda^*,Az_k-b\rangle.
	\]
	In particular, $\mathcal E_{\kappa_\tau}^{K}$
	has an upper bound independent of $\tau$.
	Applying~\eqref{AL_thm1_2} with lower summation
	index $K$ therefore yields
	\[
	\begin{aligned}
		\mathcal E_{\kappa_\tau}^{k+1}
		&\leq
		\frac{a_K}{a_{k+1}}\mathcal E_{\kappa_\tau}^{K}
		+C_1
		+\frac{\alpha\|x^*\|^2}{3a_{k+1}}
		\sum_{i=K}^{k}a_iis\epsilon_i\\
		&\leq C_2,
		\qquad k,\tau\geq K,
	\end{aligned}
	\]
	for a constant $C_2>0$ independent of $k$ and $\tau$.
	Taking $\tau=k$ gives
	\[
	(k+1)(k+2-\theta)s
	\big(
	\mathcal L(z_{k+1},\lambda^*)
	-\mathcal L(x^*,\lambda^*)
	+\|Az_{k+1}-b\|
	\big)\leq C_2.
	\]
	As in the proof of Theorem~\ref{AL_thm1},
	nonnegativity of the Lagrangian gap then implies
	\[
	\|Az_{k+1}-b\|
	\leq\frac{C_2}{(k+1)(k+2-\theta)s},
	\qquad
	|f(z_{k+1})-f(x^*)|
	\leq\frac{(1+\|\lambda^*\|)C_2}
	{(k+1)(k+2-\theta)s}.
	\]
	This completes the proof.
\end{proof}

Based on the above lemmas, we now establish the strong convergence result for the sequence generated by Algorithm \ref{Algorithm1}.

\begin{theorem}
	\label{AL_thm_strong}
	Suppose that the assumptions of Lemma~\ref{AL_lem1} hold,
	$2\leq\theta<3$, and
	$\sum_{k=K_0}^{+\infty}\epsilon_k/k<+\infty$.
	Assume that
	\[
	k^{\frac{\alpha+3}{3}}\epsilon_k\to+\infty,
	\qquad
\frac{1}{k^{\frac{\alpha+3}{3}}\epsilon_k}
\sum_{i=K_0}^{k}i^{\frac{\alpha+3}{3}}\epsilon_i^2\to0,
	\]
	and that there exists $\delta>0$ such that
	\[
	k(k+1-\theta)s\epsilon_k
	\geq\frac{(\alpha-3)(5\alpha-3)}{9}+\delta
	\]
	for all sufficiently large $k$.
	Assume additionally that $f$ is continuously differentiable
	with an $L$-Lipschitz continuous gradient.
	Let $\{(x_k,\lambda_k)\}_{k\geq K_0}$ be generated by
	Algorithm~\ref{Algorithm1}, and let
	$x^*=\operatorname{Proj}_S0$.
	Then $\liminf_{k\to+\infty}\|x_k-x^*\|=0$.
	Moreover, if $\{x_k\}$ eventually lies entirely in either
	$B(0,\|x^*\|)$ or its complement, then $x_k\to x^*$.
\end{theorem}

\begin{proof}
	By Theorem~\ref{AL_thm1} and Lemma~\ref{AL_conv_dual},
	$x_{k+1}-x_k\to0$ and $\lambda_k\to\lambda^\infty$,
	where $(x^*,\lambda^\infty)\in\Omega$.
	Here any primal solution can be paired with the limiting
	dual solution.
	Moreover, $\epsilon_k\to0$ and
	$\mathcal L_\rho(z_k,\lambda^\infty)
	-\mathcal L_\rho(x^*,\lambda^\infty)\to0$.
	We consider three cases.
	
	\medskip
	\noindent
	\textbf{Case I.}
	Suppose that $\|x_k\|\geq\|x^*\|$ for all sufficiently
	large $k$.
	
	Set $p=(\alpha-3)/3$, $\nu=p+2$, and $w_k=k^p$.
	Define
	\[
	\widetilde{\mathcal E}_k
	:=w_k\left(
	\mathcal E_{\lambda^\infty}^k
	-\frac{k(k+1-\theta)s\epsilon_k}{2}\|x^*\|^2
	\right).
	\]
	This energy is nonnegative in the present case.
	Let $d_{k+1}=x_{k+1}-x_k$ and
	$B_k=\mathcal C_k^3/(ks)$, where $\mathcal C_k^3$
	is defined in the proof of Lemma~\ref{AL_lem1}, and set
	\[
	P_k:=
	\frac{(k+1)(k+2-\theta)s}{2}(\epsilon_{k+1}-\epsilon_k)
	+\left(-pk+1-\frac{\theta}{2}\right)s\epsilon_k
	+\frac{Mk^2(k+1-\theta)|\beta_k|s^2}{4(k-\theta)}
	\epsilon_k^2.
	\]
	Taking $\kappa=\lambda^\infty$
	in~\eqref{AL_lem1_7}, and discarding the nonpositive
	penalty and squared-gradient terms for sufficiently
	large $k$, gives
	\[
	\begin{aligned}
		\widetilde{\mathcal E}_{k+1}
		-\widetilde{\mathcal E}_k
		\leq{}&
		-\mathcal A_k
		\big(
		\mathcal L_\rho(z_{k+1},\lambda^\infty)
		-\mathcal L_\rho(x^*,\lambda^\infty)
		\big)
		+\mathcal C_k
		\big(\|x_{k+1}\|^2-\|x^*\|^2\big)
		+\mathcal Q_k+R_k\|x^*\|^2\\
		&+\frac{\alpha}{3}(w_{k+1}-w_k)
		\|\lambda_{k+1}-\lambda^\infty\|^2,
	\end{aligned}
	\]
	where
	\[
	\begin{aligned}
		\mathcal A_k
		&:=2pw_k(k+2-\theta)s
		-(w_{k+1}-w_k)(k+1)(k+2-\theta)s,\\
		\mathcal C_k
		&:=w_kP_k+\frac{w_{k+1}-w_k}{2}
		(k+1)(k+2-\theta)s\epsilon_{k+1},\\
		R_k
		&:=w_k\frac{Mk^2(k+1-\theta)|\beta_k|s^2}
		{4(k-\theta)}\epsilon_k^2,
	\end{aligned}
	\]
	and
	\[
	\begin{aligned}
		\mathcal Q_k
		:={}&-\frac{\alpha}{3}w_kks\epsilon_k
		\|x_{k+1}-x^*\|^2
		+w_k\left(-p(k+1-\theta)
		+\frac{\sqrt{k}}{2}|B_k|\right)\|d_{k+1}\|^2\\
		&+(w_{k+1}-w_k)
		\left[
		\frac12\left\|\frac{2\alpha}{3}(y_{k+1}-x^*)\right\|^2
		+\frac{\alpha(\alpha-3)}{9}\|x_{k+1}-x^*\|^2
		\right].
	\end{aligned}
	\]
	
	We first show that $\mathcal Q_k\leq0$ eventually.
	Set $r_{k+1}=x_{k+1}-x^*$,
	$u_{k+1}=(k+1-\theta)d_{k+1}$,
	$\vartheta_k=(w_{k+1}-w_k)(k+1-\theta)/w_k$,
	and $\xi_k=\sqrt{k}|B_k|/[2(k+1-\theta)]$.
	Since $\{B_k\}$ is bounded,
	$\vartheta_k\to p$ and $\xi_k\to0$.
	Using
	$\frac{2\alpha}{3}(y_{k+1}-x^*)
	=\frac{2\alpha}{3}r_{k+1}+u_{k+1}$, we obtain
	\[
	\mathcal Q_k
	=\frac{w_k}{k+1-\theta}
	\left[
	-a_k\|r_{k+1}\|^2
	+\frac{2\alpha}{3}\vartheta_k
	\langle r_{k+1},u_{k+1}\rangle
	-c_k\|u_{k+1}\|^2
	\right],
	\]
	where
	\[
	a_k:=\frac{\alpha}{3}k(k+1-\theta)s\epsilon_k
	-\frac{\alpha(\alpha-1)}{3}\vartheta_k,
	\qquad
	c_k:=p-\xi_k-\frac{\vartheta_k}{2}.
	\]
	The lower bound on $k(k+1-\theta)s\epsilon_k$ implies
	\[
	\liminf_{k\to+\infty}a_k
	\geq\frac{2\alpha^2p}{9}+\frac{\alpha\delta}{3},
	\qquad
	c_k\to\frac p2,
	\]
	and hence
	\[
	\liminf_{k\to+\infty}
	\left[
	a_kc_k-\left(\frac{\alpha\vartheta_k}{3}\right)^2
	\right]
	\geq\frac{\alpha\delta p}{6}>0.
	\]
	Thus $a_k,c_k>0$ and
	$a_kc_k\geq(\alpha\vartheta_k/3)^2$ eventually.
	The quadratic form in brackets is therefore
	nonpositive, proving $\mathcal Q_k\leq0$.
	
	Next, $w_{k+1}-w_k=pk^{p-1}+o(k^{p-1})$ gives
	$\mathcal A_k=psk^{p+1}+o(k^{p+1})>0$ eventually.
	Condition~\eqref{AL_parameter_condition_1} and
	$\theta\geq2$ imply $P_k\leq-pks\epsilon_k$.
	Since $\epsilon_{k+1}\leq\epsilon_k$, it follows that
	\[
	\mathcal C_k
	\leq-pk^{p+1}s\epsilon_k
	+\frac p2 k^{p+1}s\epsilon_{k+1}
	+o(k^{p+1}\epsilon_k)
	\leq-\frac p2 k^{p+1}s\epsilon_k
	+o(k^{p+1}\epsilon_k)\leq0
	\]
	for all sufficiently large $k$.
	Using the nonnegative gap and
	$\|x_{k+1}\|^2-\|x^*\|^2\geq0$, we conclude that
	\[
	\widetilde{\mathcal E}_{k+1}
	-\widetilde{\mathcal E}_k
	\leq R_k\|x^*\|^2
	+\frac{\alpha}{3}(w_{k+1}-w_k)
	\|\lambda_{k+1}-\lambda^\infty\|^2.
	\]
	Because $\{\beta_k\}$ is bounded,
	$R_k\leq Ck^\nu\epsilon_k^2$ and
	$w_{k+1}-w_k\leq Ck^{\nu-3}$ eventually.
	Summing from a sufficiently large integer $K_1$ to $k-1$
	gives
	\[
	\widetilde{\mathcal E}_k
	\leq\widetilde{\mathcal E}_{K_1}
	+C\sum_{i=K_1}^{k-1}i^\nu\epsilon_i^2
	+C\sum_{i=K_1}^{k-1}i^{\nu-3}
	\|\lambda_{i+1}-\lambda^\infty\|^2.
	\]
	
	Let $x_{\epsilon_k}$ be the regularized minimizer
	defined in Lemma~\ref{lem:tikhonov_AL}.
	Then $x_{\epsilon_k}\to x^*$.
	By strong convexity of
	$\mathcal L_\rho(\cdot,\lambda^\infty)
	+\frac{\epsilon_k}{2}\|\cdot\|^2$ and
	$z_k=x_k+\beta_k(x_k-x_{k-1})$,
	\[
	\begin{aligned}
		\widetilde{\mathcal E}_k
		\geq{}&
		\frac{w_kk(k+1-\theta)s\epsilon_k}{2}
		\Big(
		\|x_k-x_{\epsilon_k}\|^2
		+\|x_{\epsilon_k}\|^2-\|x^*\|^2
		-2\beta_k\langle x_k-x_{k-1},x_{\epsilon_k}\rangle
		\Big).
	\end{aligned}
	\]
	Combining these bounds and using
	$w_kk(k+1-\theta)s\sim sk^\nu$, we obtain
	\[
	\begin{aligned}
		\|x_k-x_{\epsilon_k}\|^2
		\leq{}&
		\|x^*\|^2-\|x_{\epsilon_k}\|^2
		+2\beta_k\langle x_k-x_{k-1},x_{\epsilon_k}\rangle\\
		&+\frac{C}{k^\nu\epsilon_k}
		+\frac{C}{k^\nu\epsilon_k}
		\sum_{i=K_1}^{k-1}i^\nu\epsilon_i^2
		+\frac{C}{k^\nu\epsilon_k}
		\sum_{i=K_1}^{k-1}i^{\nu-3}
		\|\lambda_{i+1}-\lambda^\infty\|^2.
	\end{aligned}
	\]
	The first two terms tend to zero because
	$x_{\epsilon_k}\to x^*$, $x_k-x_{k-1}\to0$,
	and $\{\beta_k\}$ is bounded.
	The next two terms vanish by the assumptions on
	$\epsilon_k$.
	For the last term, the lower bound on
	$k(k+1-\theta)s\epsilon_k$ gives
	\[
	0\leq
	\frac{\displaystyle\sum_{i=K_1}^{k-1}
		i^{\nu-3}\|\lambda_{i+1}-\lambda^\infty\|^2}
	{k^\nu\epsilon_k}
	\leq
	\frac{C}{k^{\nu-2}}
	\sum_{i=K_1}^{k-1}
	i^{\nu-3}\|\lambda_{i+1}-\lambda^\infty\|^2
	\to0.
	\]
	The last limit follows by splitting the sum into a finite
	initial part and a tail, since $\nu>2$ and
	$\lambda_k\to\lambda^\infty$.
	Thus $\|x_k-x_{\epsilon_k}\|\to0$, and consequently
	$x_k\to x^*$.
	
	\medskip
	\noindent
	\textbf{Case II.}
	Suppose that $\|x_k\|<\|x^*\|$ for all sufficiently
	large $k$.
	
	Then $\{x_k\}$ is bounded.
	Since $z_k-x_k=\beta_k(x_k-x_{k-1})\to0$,
	continuity of $\mathcal L_\rho(\cdot,\lambda^\infty)$
	on bounded sets and the gap convergence imply
	\begin{equation}\label{AL_thm3.3_ex1}
		\mathcal L_\rho(x_k,\lambda^\infty)
		-\mathcal L_\rho(x^*,\lambda^\infty)\to0.
	\end{equation}
	Let $\widetilde x$ be any cluster point.
	Passing to the limit along a convergent subsequence gives
	$\mathcal L_\rho(\widetilde x,\lambda^\infty)
	=\mathcal L_\rho(x^*,\lambda^\infty)$ and
	$\|\widetilde x\|\leq\|x^*\|$.
	As shown in Lemma~\ref{lem:tikhonov_AL},
	$\operatorname*{argmin}\mathcal L_\rho
	(\cdot,\lambda^\infty)=S$.
	Hence $\widetilde x\in S$, and uniqueness of the
	minimum-norm solution gives $\widetilde x=x^*$.
	The bounded sequence therefore converges to $x^*$.
	
	\medskip
	\noindent
	\textbf{Case III.}
	Suppose that neither Case I nor Case II holds.
	
	There exists a subsequence $\{x_{k_j}\}$ with
	$\|x_{k_j}\|<\|x^*\|$.
	It is bounded and admits a convergent further subsequence.
	By the argument of Case II, its limit belongs to $S$
	and has norm at most $\|x^*\|$, so the limit is $x^*$.
	Therefore $\liminf_{k\to+\infty}\|x_k-x^*\|=0$.
	This completes the proof.
\end{proof}

\section{Numerical Experiments}\label{Numerical}
In this section, we report numerical experiments to validate the theoretical results developed in the preceding sections. All experiments are implemented in \textsc{Matlab} (R2021b) and executed on a PC equipped with a 12th Gen Intel(R) Core(TM) i7-12700 CPU at 2.10~GHz and 32~GB of RAM.

\noindent\textbf{Example 1: Convex optimization problem.}
\begin{equation*}
	\begin{aligned}
		&\min_{x}\,(mx_1+nx_2+ex_3)^2,\\
		&\;s.t. \, Ax=b,
	\end{aligned}
\end{equation*}
where $A=(m, -n, e)$ with $m, n, e\in\mathbb{R}$\textbackslash$\left\{0\right\}$ and $b=0$. A straightforward calculation yields that the solution set can be parameterized as $S=\left\{\, (x_1,0,-\frac{m}{e}x_1)^{\top}\ \big|\ x_1\in\mathbb{R}\,\right\}$. It follows that the unique solution of minimum-norm is $x^{\ast}=(0,0,0)^{\top}$, and the corresponding optimal objective value equals \(0\).

\emph{Experiment 1: Sensitivity with respect to the parameters $\gamma$ and $\beta$.} 
For the mixed primal--dual dynamical system with implicit Hessian damping and Tikhonov regularization \eqref{DS} (abbreviated as MPDDS-IHDTR), we investigate the influence of the parameters \(\gamma\) and \(\beta\) on the numerical performance. Specifically, by varying \(\gamma\) and \(\beta\), we examine how the \emph{feasibility violation} and the \emph{primal--dual gap} change with respect to these parameters. The parameter settings used in the experiments are summarized as follows:

\noindent $\blacktriangleright$ Dynamical system (MPDDS-IHDTR). We fix
$ m=1, n=10^{-3}, e=10^{-4}, \alpha=3.1, \rho=1, \epsilon(t)=t^{-r}, r=1.1$, and consider the following pairs \((\gamma,\beta)\):
$ (\gamma,\beta)\in\{(0,0),\ (0,1),\ (1,-0.5),\ (1,0.5)\}$.

\begin{figure}[H]
	\centering
	{
		\begin{minipage}[t]{0.48\linewidth}
			\centering
			\includegraphics[width=2.65in]{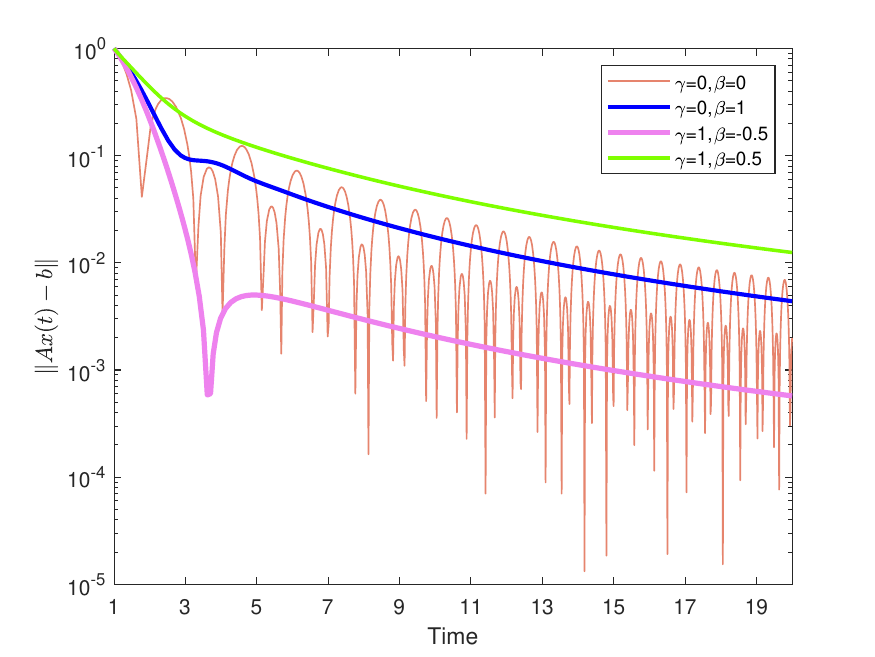}
		\end{minipage}%
	}
	{
		\begin{minipage}[t]{0.48\linewidth}
			\centering
			\includegraphics[width=2.65in]{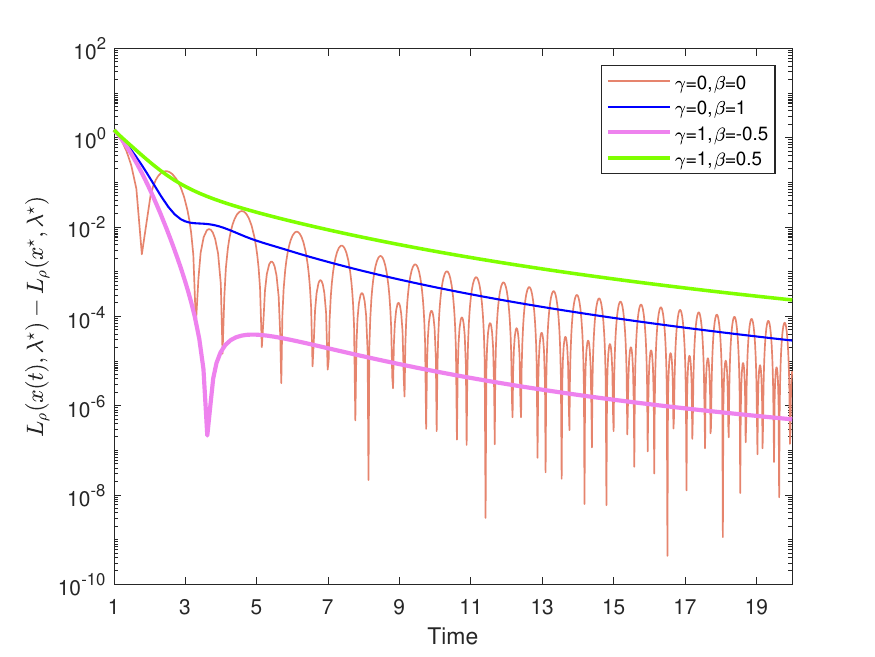}
		\end{minipage}%
	}
	\caption{Error analysis with different values of parameters $\gamma$ and $\beta$ in the dynamical system \eqref{DS} for Example 1.}
	\label{fig: 1_self}
	\centering
\end{figure}

\vspace{-0.8cm}
Figure \ref{fig: 1_self} reports the time evolution of the feasibility violation and the primal--dual gap for different choices of \((\gamma,\beta)\).
Overall, all tested configurations exhibit a clear decay of both criteria, indicating convergence of the trajectories toward a primal--dual solution.
However, the convergence speed and stability are markedly sensitive to \((\gamma,\beta)\). In particular, the baseline setting \((\gamma,\beta)=(0,0)\) shows pronounced oscillations, whereas increasing \(\beta\) (e.g., \((0,1)\)) noticeably reduces oscillatory behavior and improves the decay profile. When \(\gamma\) is activated, the choice of \(\beta\) becomes crucial: a negative \(\beta\) yields a substantially faster reduction of both feasibility violation and primal--dual gap, while a positive \(\beta\) leads to a slower decay. These observations suggest that incorporating implicit Hessian damping together with a suitable Tikhonov regularization coupling improves numerical stability and accelerates convergence, in agreement with the theoretical insights established earlier.

\emph{Experiment 2: Comparison among five dynamical systems on Example~1.} 
We further compare the numerical performance of five dynamical systems for solving Example~1: (i) the mixed primal--dual dynamical system \eqref{DS} with implicit Hessian damping and Tikhonov regularization (MPDDS-IHDTR), (ii) the system \eqref{DS} without the damping term $\bigl(\gamma+\tfrac{\beta}{t}\bigr)\dot{x}(t)$ (MPDDS-TR), (iii) the system \eqref{DS} without the Tikhonov regularization term $\epsilon(t)x(t)$ (MPDDS-IHD), (iv) the system \eqref{DS} without both $\bigl(\gamma+\tfrac{\beta}{t}\bigr)\dot{x}(t)$ and $\epsilon(t)x(t)$ (MPDDS), and (v) the primal--dual accelerated vanishing damping dynamics (PD-AVD) proposed in~\cite{Bot JDE 2021}.
The parameter settings are summarized as follows.
\begin{itemize}
	\item Dynamical system (MPDDS): $m=1, n=10^{-3}, e=10^{-4}, \alpha=3.1 \text{ and }\rho=1$.
	\item Dynamical system  (MPDDS-IHDTR): the same values of $(m,n,e,\alpha,\rho)$ as in the dynamical system (MPDDS), $\gamma=1,\beta=-0.5 \text{ and }\epsilon(t)=t^{-r} \text{ with }r=1.1$.
	\item Dynamical system (MPDDS-TR): the same values of $(m,n,e,\alpha,\rho)$ as in the dynamical system (MPDDS), and $ 
	\epsilon(t)=t^{-r }\text{ with }r=1.1$.
	\item Dynamical system (MPDDS-IHD): the same values of $(m,n,e,\alpha,\rho)$ as in the dynamical system (MPDDS), $\gamma=1 \text{ and } \beta=-0.5$.
	\item Dynamical system (PD-AVD): the same values of $(m,n,e,\alpha)$ as in the dynamical system (MPDDS), $\beta=1 \text{ and }\delta=\frac{1}{2.05}$.
\end{itemize}

\begin{figure}[H]
	\centering
	{
		\begin{minipage}[t]{0.32\linewidth}
			\centering
			\includegraphics[width=2.15in]{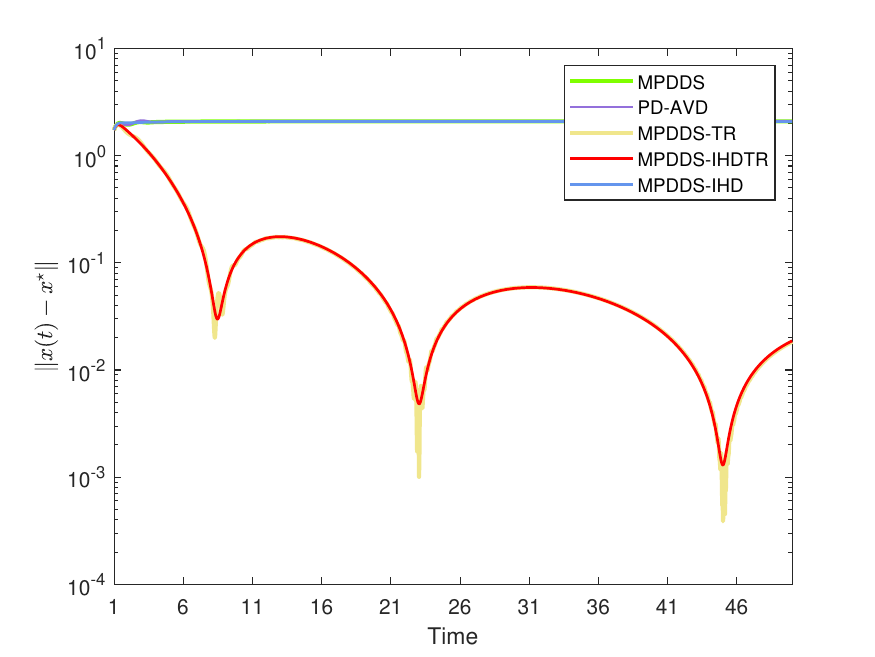}
		\end{minipage}%
	}
	{
		\begin{minipage}[t]{0.32\linewidth}
			\centering
			\includegraphics[width=2.15in]{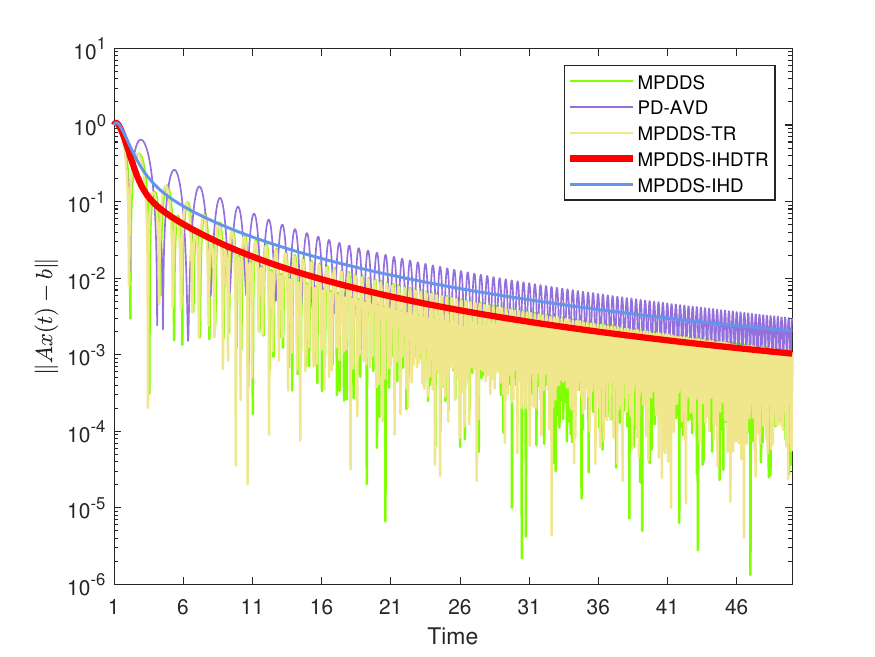}
		\end{minipage}%
	}
	{
		\begin{minipage}[t]{0.32\linewidth}
			\centering
			\includegraphics[width=2.15in]{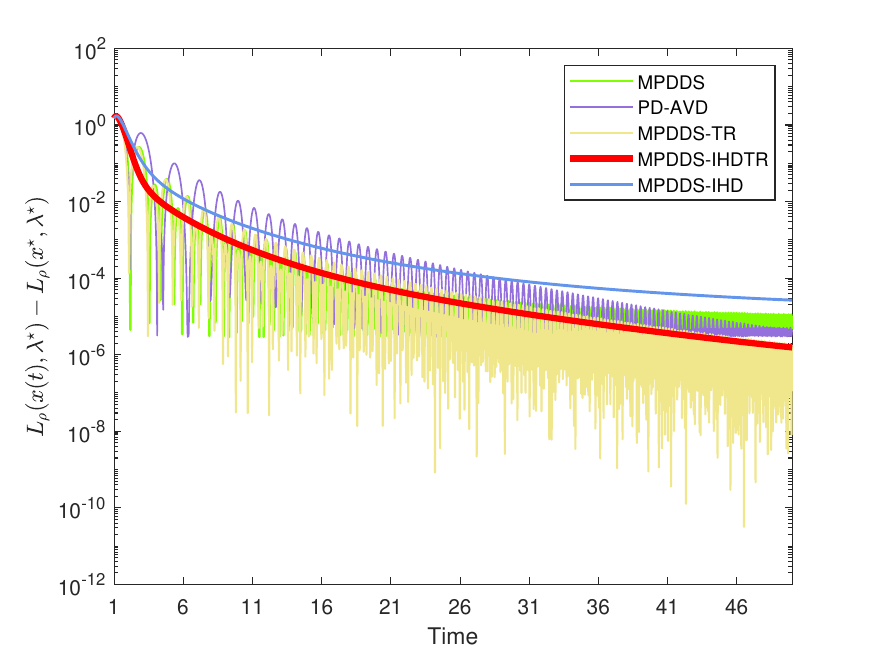}
		\end{minipage}%
	}
	\caption{Comparison of error results between (PD-AVD), (MPDDS-IHDTR), (MPDDS-IHD), (MPDDS-TR) and  (MPDDS).}
	\label{fig: 2_comparasion}
	\centering
\end{figure}

As shown in Figure~\ref{fig: 2_comparasion}, the first subplot
shows that the two Tikhonov-regularized dynamical systems,
\textup{(MPDDS-IHDTR)} and \textup{(MPDDS-TR)}, attain much
smaller errors than the other three dynamical systems.
This suggests that Tikhonov regularization promotes convergence
toward the minimum-norm solution, consistent with its role
in our strong convergence analysis.
The second and third subplots show pronounced oscillations
in the feasibility violation and primal--dual gap for
\textup{(PD-AVD)}, \textup{(MPDDS-TR)}, and \textup{(MPDDS)},
which lack implicit Hessian damping.
Both quantities decay more smoothly for
\textup{(MPDDS-IHDTR)} and \textup{(MPDDS-IHD)},
illustrating the damping effect.
Finally, \textup{(MPDDS-IHDTR)} achieves the best final accuracy
among the five dynamical systems in this experiment,
suggesting that combining Tikhonov regularization with
implicit Hessian damping improves numerical performance.

\vskip 5pt
\noindent\textbf{Example 2: Linearly constrained $\mathbf{\ell_1}$--$\mathbf{\ell_2}$ optimization problem.}
\begin{equation*}
	\begin{aligned}
		&\min_{x}\;\Vert x\Vert_1+\frac{\delta}{2}\Vert x\Vert_2,\\
		&\,s.t.\; Ax=b,
	\end{aligned}
\end{equation*}
where $A\in\mathbb{R}^{m\times n}, b\in \mathbb{R}^n$. 

We set $m=1850$, $n=180$, and $\delta=1$. The sensing matrix $A\in\mathbb{R}^{m\times n}$ is generated with entries that are independent and identically distributed according to the standard Gaussian distribution. The original solution (signal) $x^\ast\in\mathbb{R}^n$ is constructed to be sparse and contains $5\%$ nonzero components. The nonzero entries are sampled from $\mathcal{N}(0,4)$ and are then truncated to the interval $[-2,2]$. The noise vector $\omega$ is sampled from the standard Gaussian distribution and scaled to satisfy $\|\omega\|=10^{-6}$. The observation vector is given by
\begin{equation*}
	b = A x^\ast + \omega .
\end{equation*}

To solve the subproblems, we apply FISTA~\cite{Beck SIAM JIS 2009}. The inner iterations are terminated once
\begin{equation*}
	\frac{\Vert z_k-z_{k-1}\Vert}{\max \left\{\Vert z_{k-1}\Vert, 1\right\}}\leq 10^{-10}
\end{equation*}
or when the number of iterations reaches $100$, where $\{z^k\}$ denotes the iterative sequence produced by FISTA. 

Next, we compare Algorithm \ref{Algorithm1}, the accelerated linearized augmented Lagrangian method (ALALM)~\cite[Algorithm 1 with adaptive parameters]{Xu SIOPT 2017}, and the inertial accelerated primal--dual method (IAPDM)~\cite[Algorithm~1]{He NA 2022} for solving Example~2. In particular, we examine how the relative error ($\mathrm{Rel}$) and the residual error ($\mathrm{Res}$) evolve during the iterations. The parameter settings of these three algorithms are specified as follows:

\begin{itemize}
	\item Algorithm \ref{Algorithm1}: $\alpha=9, \theta=2.1, \rho=0.1, \gamma=1, \beta=-0.5, s=1, \epsilon_k=\frac{1}{k^3}$.\\
	\item ALALM: $\gamma =1, \alpha_k=\frac{2}{k+1}, \gamma_k=k\gamma, \beta_k=\frac{\gamma_k+1}{2}$ and $p_k=\frac{1}{k}\rm{Id}$.\\
	\item IAPDM: $\alpha =9, s=1, M_k=\frac{1}{m}Id$.
\end{itemize}
\begin{figure}[H]
	\centering
	{
		\begin{minipage}[t]{0.48\linewidth}
			\centering
			\includegraphics[width=2.56in]{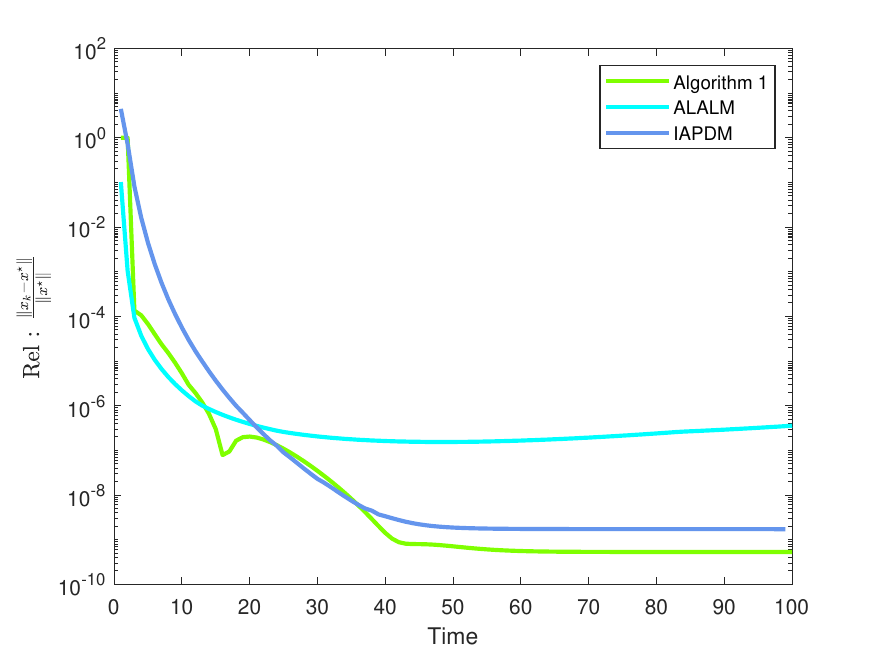}
		\end{minipage}%
	}
	{
		\begin{minipage}[t]{0.48\linewidth}
			\centering
			\includegraphics[width=2.56in]{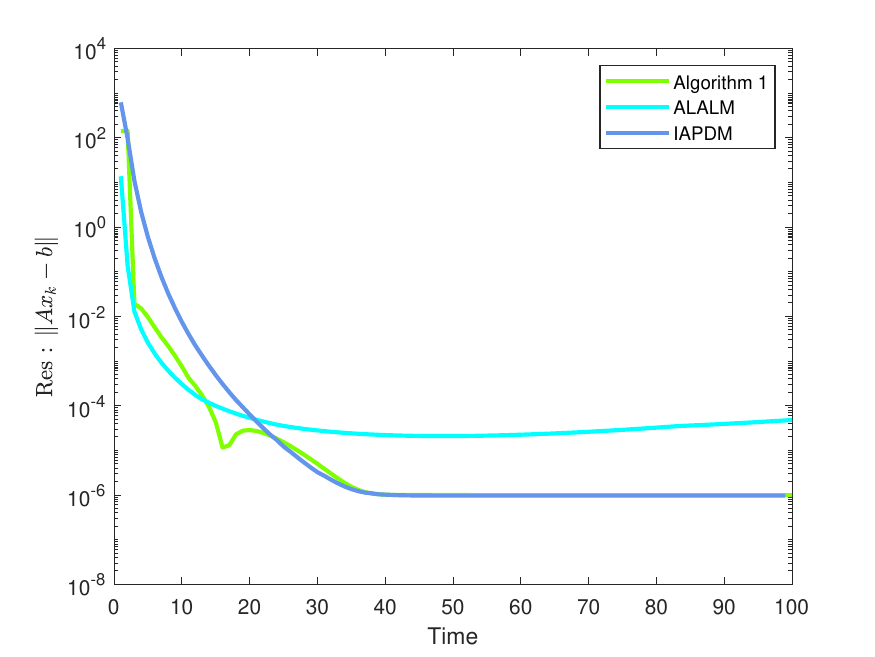}
		\end{minipage}%
	}
	\caption{Comparisons of Algorithm \ref{Algorithm1}, ALALM and IAPDM for Example 2.}
	\label{fig: 3_algorithm}
	\centering
\end{figure}

From Figure \ref{fig: 3_algorithm}, we observe that Algorithm \ref{Algorithm1} exhibits a faster decrease overall. In particular, in the first subplot of Figure \ref{fig: 3_algorithm} (the evolution of $\mathrm{Rel}$), Algorithm \ref{Algorithm1} decreases at a rate comparable to ALALM during the early iterations, while it attains a higher solution accuracy in the later stage. Compared with IAPDM, Algorithm \ref{Algorithm1} is clearly superior in terms of $\mathrm{Rel}$, achieving both a faster reduction and a smaller final error level. We next consider the second subplot of Figure \ref{fig: 3_algorithm}, which reports the evolution of $\mathrm{Res}$. It can be seen that, relative to ALALM, Algorithm \ref{Algorithm1} still has the advantage of reaching a lower residual level, indicating better final accuracy. Moreover, when compared with IAPDM, Algorithm \ref{Algorithm1} shows a noticeably faster decrease of $\mathrm{Res}$. Taken together, these results indicate that Algorithm \ref{Algorithm1} not only improves the attainable solution accuracy but also achieves this accuracy with fewer iterations, thereby leading to a faster overall solution process.

\section{Conclusion}\label{Conclusion}

This paper proposed a class of mixed primal--dual dynamical systems with implicit Hessian damping and Tikhonov regularization for solving problem~\eqref{problem1}. By constructing suitable energy functions, we established accelerated $\mathcal{O}(1/t^{2})$ convergence rates for the primal--dual gap, the feasibility violation, and the objective residual along
$x(t)+\big(\gamma+\frac{\beta}{t}\big)\dot{x}(t)$, and we proved that the generated trajectory converges strongly to the minimum-norm solution of~\eqref{problem1}. 
The proposed model effectively reduces oscillations without requiring the objective function to be twice differentiable. Moreover, although the dual variable follows a first-order ordinary differential equation, the overall system still achieves an accelerated rate, and its mixed ``second-order'' plus ``first-order'' structure can be attractive in terms of computational cost. We further derived a fast primal--dual algorithm via an implicit discretization of~\eqref{DS} and showed, through discrete energy sequences, that the algorithm inherits the accelerated rate and strong convergence. Numerical experiments validated the theoretical findings and demonstrated fast decrease and high solution accuracy in practice.

\vspace{5mm}
\noindent \textbf{Acknowledgements} This research was supported by the National Natural Science Foundation of China (12571318, 12171063) and Sichuan Science and Technology Program (Grant No. 2025ZNSFSC0813).

\section*{Declarations}
\textbf{Competing interest} The authors declare that they have no competing interests.

\end{document}